\documentclass[smallextended]{svjour3}       
\smartqed  
\usepackage{graphicx}

\usepackage{amsmath}
\usepackage{enumerate}

\usepackage{amssymb}
\usepackage{graphicx}
\usepackage{epstopdf}
\usepackage{subfigure}
\usepackage{color}
\usepackage{array}
\usepackage{mathrsfs}
\usepackage{multirow}
\usepackage{algorithm,algorithmic}
\usepackage{booktabs}
\usepackage{multicol}

\begin{document}

\title{A high-order discrete energy decay and maximum-principle preserving scheme for time fractional Allen-Cahn equation
}

\titlerunning{A high-order scheme for time fractional Allen-Cahn equation}
\author{Guoyu Zhang$^1$      \and
        Chengming Huang$^{2,3}$   \and
        Anatoly A. Alikhanov$^4$     \and
        Baoli Yin$^1$
}


\institute{G. Zhang  \\
              \email{guoyu\_zhang@imu.edu.cn}  \\
              C. Huang \\
              \email{chengming\_huang@hotmail.com} \\         
              A. Alikhanov\\
               \email{aaalikhanov@gmail.com}\\
             B. Yin\\
               \email{baolimath@126.com}
           \at
           $^1$School of Mathematical Sciences, Inner Mongolia University, Hohhot 010021, China \\
           $^2$School of Mathematics and Statistics, Huazhong University of Science and Technology, Wuhan 430074, China \\
           $^3$Hubei Key Laboratory of Engineering Modeling and Scientific Computing, Huazhong University of Science and Technology, Wuhan 430074, China \\
           $^4$North-Caucasus Center for Mathematical Research, North-Caucasus Federal University, Stavropol 355017, Russia
}

\date{Received: date / Accepted: date}

\maketitle

\begin{abstract}
The shifted fractional trapezoidal rule (SFTR) with a special shift is adopted to construct a finite difference scheme for the time-fractional Allen-Cahn (tFAC) equation. Some essential key properties of the weights of SFTR are explored for the first time. Based on these properties, we rigorously demonstrate the discrete energy decay property and maximum-principle preservation for the scheme. Numerical investigations show that the scheme can resolve the intrinsic initial singularity of such nonlinear fractional equations as tFAC equation on uniform meshes without any correction. Comparison with the classic fractional BDF2  and L2-1$_\sigma$ method further validates the superiority of SFTR in solving the tFAC equation. Experiments concerning both discrete energy decay and discrete maximum-principle also verify the correctness of the theoretical results.

\end{abstract}

\section{Introduction}\label{intro}
For systems driven by dissipation of free energy, the gradient flow is a useful tool in diverse research areas (see, e.g., \cite{AndersonMcFaddenWheeler,AllenCahn,ShaoRappelLevine,XuPrinceSnakes,LiuChengWang}). 
In this work, we are interested in systems where dissipation of the free energy has memory effect and thus can be simulated by the $\alpha$th order ($0<\alpha<1$) fractional gradient flow (see \cite{HouZhuXu,DuYangZhou,TangYuZhou}),
\begin{equation}\label{Int.0}\begin{split}
\partial_t^\alpha u(\boldsymbol x,t)=-\frac{\delta E}{\delta u},\quad
\boldsymbol x\in \Omega,~t\in (0,T],
\end{split}
\end{equation}
with the energy 
\begin{equation*}\begin{split}
E[u](t)=\int_{\Omega}\bigg(\frac{\varepsilon^2}{2}|\nabla u|^2+F(u)\bigg)\mathrm{d}\boldsymbol x,
\end{split}
\end{equation*}
where $T>0$, $\Omega=(a,b)\times(a,b)$, $F(u)=\frac{1}{4}(1-u^2)^2$ is a double-well potential, and the constant $\varepsilon(0<\varepsilon\ll1)$ represents the thickness of the transition boundary between materials.
The symbol $\partial_t^\alpha$ denotes the Caputo fractional operator known as
\begin{equation}\label{Int.0.1}\begin{split}
\partial_t^\alpha u(t)=\frac{1}{\Gamma(1-\alpha)}\int_{0}^{t}\frac{u'(s)}{(t-s)^\alpha}\mathrm{d}s,
\quad
\Gamma(\alpha) \text{ is the Gamma function}.
\end{split}
\end{equation}
Imposed the periodic boundary condition, equation (\ref{Int.0}) is then essentially the nonlocal phase field model, or specifically, the time-fractional Allen-Cahn (tFAC) equation
\begin{equation}\label{Int.1}\begin{split}
\partial_t^\alpha u(\boldsymbol x,t)&=\varepsilon^2\Delta u-f(u),\quad
\boldsymbol x\in \Omega,~
t\in (0,T],
\quad
u(\boldsymbol x,0)=u_0(\boldsymbol x),~
\boldsymbol x\in \overline{\Omega},
\end{split}
\end{equation}
where $f(u)=F'(u)=u^3-u$ represents  the nonlinear bulk force.
The nonlocality of the fractional derivative and nonlinearity of equation (\ref{Int.1}) generally render solving this problem rather difficult by analytic methods.
We thus consider numerical methods for this problem. In this context, an important issue is to investigate whether numerical methods can inherit some intrinsic properties of the solution, including the energy decay and maximum-principle property. 
\par
It is well known that for classical Allen-Cahn equation (with $\alpha=1$), the energy $E[u](t)$ admits a monotonic property:
\begin{equation}\label{Int.1.1}\begin{split}
E'[u](t)+\bigg\|\frac{\delta E}{\delta u}\bigg\|^2=0,\quad \text{or}\quad
E[u](t_1)\leq E[u](t_2), \quad\text{if } t_1\geq t_2,
\end{split}
\end{equation}
while for general cases (with $\alpha \in (0,1)$), the monotonicity (\ref{Int.1.1}) still remains an open problem, despite the fact that all numerical tests published in the literature up to now suggest its correctness.
To circumvent this difficulty, on the one hand a ``weak version'' of $(\ref{Int.1.1})$ called energy dissipation, i.e. $E[u](t)\leq E[u](0)$ for $t \geq 0$, was studied in both continuous and discrete levels by Tang et al.\cite{TangYuZhou} for time-fractional phase-field models (e.g., the time-fractional Allen-Cahn/Cahn-Hilliard/molecular-beam-epitaxy equation).
See also \cite{ChenZhangZhaoCaoWangZhang}.
On the other hand, some modified types of energy have been proposed for the model (\ref{Int.1}), in which the variable $\alpha$ is involved.
Specifically, the authors in \cite{QuanTangYang} defined a weighted energy
\begin{equation}\label{Int.1.2}\begin{split}
E_w[u](t)=\int_{0}^{1}w(s)E[u](st)\mathrm{d}s,
\end{split}
\end{equation}
where $w(s)$ denotes some weight function satisfying $\int_{0}^{1}w(s)\mathrm{d}s=1$, and proved $E'_w[u](t)\leq 0$ provided $w(s)s^{1-\alpha}(1-s)^\alpha$ is nonincreasing w.r.t. $s$.
It is notable that the energy $E_w[u](t)$ involves $\alpha$ since $w(s)$ in general depends on $\alpha$.
Considering the tFAC equation degrades to the classical case, it is reasonable to require the modified energy recover some aspects of (\ref{Int.1.1}) when $\alpha \to 1$, leading to the compatible energy \cite{LiaoTangZhou1,LiaoZhuWang},
\begin{equation}\label{Int.1.3}\begin{split}
E_c[u]=E[u]+\frac{1}{2}\mathcal{I}^\alpha_t\bigg\|\frac{\delta E}{\delta u}\bigg\|^2,
\end{split}
\end{equation}
where $\mathcal{I}^\alpha_t$ denotes the Riemann-Liouville fractional integral operator of order $\alpha$ defined by $(\mathcal{I}^\alpha_t u)(t)=(\kappa*u)(t)$ with the kernel function $\kappa(s)=s^{\alpha-1}/\Gamma(\alpha)$ and $*$ standing for the convolution.
With this modified energy $E_c[u]$, one can get
\begin{equation}\label{Int.1.4}\begin{split}
E'_c[u](t)+\frac{t^{\alpha-1}}{2\Gamma(\alpha)}\bigg\|\frac{\delta E}{\delta u}\bigg\|^2 \leq 0 \quad\text{for } t>0,
\end{split}
\end{equation}
which degrades into the form $E'[u](t)+\big\|\frac{\delta E}{\delta u}\big\|^2 \leq 0$ if $\alpha \to 1$, c.f. (\ref{Int.1.1}) for the classical Allen-Cahn equation. Other types of modified energy can be found, for example, in \cite{hou2022second,hou2021highly}. 
In addition to energy decay, another important property of  (\ref{Int.1}) is that its solution satisfies the maximum-principle \cite{TangYuZhou}, i.e.,
\begin{equation}\label{Int.1.5}\begin{split}
|u(\boldsymbol x,t)|\leq 1,\quad \text{if }~|u_0(\boldsymbol x)|\leq 1.
\end{split}
\end{equation}
Note that the classicial model ($\alpha=1$) also holds such property as (\ref{Int.1.5}) \cite{ShenTangYang}.
\par
In the discrete case, it is hence necessary to develop numerical schemes that preserve the energy decay and maximum-principle for model (\ref{Int.1}).
Historically, there has been enormous numerical exploration for tFAC equations or phase field models or more generally, the nonlinear subdiffustion models, where some energy stable schemes were developed (see, e.g.,  \cite{JiaZhangXuJiang,ji2020adaptive,huang2020optimal,li2021conforming,li2019nonconforming,jin2018numerical,bu2017finite,li2018finite,liu2009numerical,zhang2020efficient,mohebbi2013high,li2018unconditionally,gao2011compact}).
However, the studies of preservation on discrete energy dissipation, or energy decay and maximum-principle are still limited.
Liao et al.\cite{LiaoTangZhou} presented a second-order nonuniform time-stepping scheme for (\ref{Int.1}) and demonstrated the scheme can preserve the discrete maximum-principle.
In \cite{MaskariKaraa}, the authors presented a finite element method for the time fractional Cahn-Hilliard equation using the fractional Euler formula in temporal direction and proved the scheme can preserve the energy dissipation. Very recently, the authors in \cite{LiaoTangZhou1} constructed a $(1+\alpha)$th order Crank-Nicolson scheme for the tFAC equation by applying an L1-type formula of the Riemann-Liouville derivative, which can preserve both the discrete version energy decay (\ref{Int.1.4}) and maximum-principle (\ref{Int.1.5}).
Hou and Xu \cite{hou2022second} explored the time
fractional $L^2$ gradient flows where the fractional derivative was approximated by the well known L2-1$_\sigma$ method \cite{alikhanov2015new}  coupled with the so-called sum-of-exponentials (SOE) technique \cite{jiang2017fast} and proved a modified type of energy dissipation  on uniform meshes.
They obtained the optimal second-order accuracy numerically by using graded meshes to compensate the low accuracy at initial time.
However, to the best of our knowledge, high-order methods on uniform meshes that can preserve simultaneously the two intrinsic properties, i.e., the energy decay and maximum-principle, of tFAC equation are still missing in the literature.
In this work, we shall fill this gap by a closer examination of a second-order accuracy formula called shifted fractional trapezoidal rule (SFTR) \cite{Yin01,YinLiuLiZhang222}, revealing some essential characteristics of this formula that are not shared with for most other second-order ones.
\par
It has been a consensus that the solutions of time fractional models in general have intrinsic initial singularities, which may lead to order reduction for numerical methods of high order designed for smooth solutions 
\cite{Stynes,zhao2014collocation}.
To resolve this issue among finite difference methods, two popular techniques are often considered, i.e., adding correction terms if the mesh is uniformly divided (see e.g. \cite{JinLazarovZhou,Lubich1}), or adopting nonuniform meshes that are much denser near the initial point (see e.g. \cite{LiaoTangZhou,stynes2017error,kopteva2019error}, to name just a few).
Some other techniques, such as constructing log type basis functions \cite{chen2020spectrally}, or employing a nonlinear transformation for the time variable before discretization \cite{li2021novel}, are also developed recently.
Although our interest in this work mainly focuses on the study of the preservation of energy decay and maximum-principle, 
we demonstrate in numerical tests that our scheme, without any correction term, can result optimal accuracy at any positive time, even on uniform meshes.
We argue that this is the first report of such phenomenon, i.e., optimal accuracy can be automatically recovered,  for nonlinear time-fractional models.
For a similar phenomenon examined for linear models, readers can refer \cite{YinLiuLiZhang222}.
\par
We sum up the main contributions of this work:
\begin{itemize}
  \item Prove theoretically some essential properties of the weights of SFTR-$\frac{1}{2}$ (i.e., SFTR with the shift parameter as $\frac{1}{2}$) and further demonstrate the discrete energy decay and maximum-principle preservation of the scheme.
  \item Demonstrate numerically the superiority of SFTR-$\frac{1}{2}$ that the intrinsic intial singularity can be resolved  automatically.
\end{itemize}
\par
The outline of the rest of the paper is as follows.
In section \ref{sec.Novel}, we prove the key properties of the weights of SFTR-$\frac{1}{2}$ which are crucial for the discrete maximum-principle preservation.
In section \ref{sec.Energymax}, we first present the Crank-Nicolson finite difference scheme for the tFAC equation and then introduce a sequence closely related the initial weights of SFTR-$\frac{1}{2}$ with rigorous justification of its properties being essential to the discrete energy decay.
At last, we prove the discrete maximum-princple.
In section \ref{sec.tests}, we first examine our theoretical results, the discrete energy decay and maximum-principle preservation, by resorting to simulating the model (\ref{Int.1}) with a randomly distributed initial condition, and then verify the high accuracy of SFTR-$\frac{1}{2}$ formula, by comparisons with the classical fractional BDF2 and L2-1$_\sigma$ method, that our method can automatically resolve the initial singularity.
Some concluding remarks are made in section \ref{sec.conc}.

\section{SFTR and kernel properties}\label{sec.Novel}
Divide the time interval by $0=t_0<t_1<t_2<\cdots<t_N=T$ with $t_n=n\tau$, $\tau=T/N$.
Set $\phi^n=\phi(t_n)$.
At the node $t_{n-\frac{1}{2}}$, the SFTR-$\frac{1}{2}$ formula for the fractional derivative $\partial_t^\alpha \phi$ takes the form (see \cite{Yin01,YinLiuLiZhang222})
\begin{equation}\label{key.1}\begin{split}
\partial_\tau^{\alpha,n-\frac{1}{2}} \phi=\tau^{-\alpha}\sum_{m=0}^{n}\omega_m (\phi^{n-m}-\phi^0),
\end{split}
\end{equation}
where the weights $\omega_m$ can either be obtained by the generating function
\[
\omega(\xi)=\sum_{m=0}^{\infty}\omega_m \xi^m=\bigg[\frac{1-\xi}{\frac{1}{2}(1+\xi)+\frac{1}{2\alpha}(1-\xi)}\bigg]^\alpha,
\]
or calculated directly by the recursive relation
\begin{equation}\label{key.3}\begin{split}
\omega_0&=\bigg(\frac{2\alpha}{\alpha+1}\bigg)^{\alpha},\quad
\omega_1=-\alpha\bigg(\frac{2\alpha}{\alpha+1}\bigg)^{\alpha+1},
\\
\omega_m&=\frac{2\alpha}{m(\alpha+1)}\bigg\{\bigg[\frac{1}{\alpha}(m-1)-\alpha\bigg]\omega_{m-1}+\frac{\alpha-1}{2\alpha}(m-2)\omega_{m-2}\bigg\}, m\geq 2.
\end{split}
\end{equation}
It has been shown in \cite{Yin01} the formula (\ref{key.1}) approximates $\partial_t^\alpha\phi(t_{n-\frac{1}{2}})$ with second-order accuracy. In fact,  it is a special case of general SFTR introduced in \cite{Yin01,YinLiuLiZhang222}. In the following, we give some peculiar characteristics of this formula.
\begin{lemma}\label{lem.2.1}
The following defined function $\Phi(x,y)$ is nonnegative for $(x,y)\in [0,1]\times[3,+\infty)${\rm:}
\begin{equation}\label{key.2.1}\begin{split}
\Phi(x,y)=\frac{1}{y}(y-1-x^2)\bigg(1-\frac{x+1}{2}\frac{y-1}{y}\bigg)-\frac{(1-x)(y-1)}{2(y-x^2)}\bigg(1+\frac{1}{y}\bigg).
\end{split}
\end{equation}
\end{lemma}
\begin{proof}
Replace $y$ with $3/\widetilde{y}$ for $\widetilde{y}\in(0,1]$ and set $\widetilde{\Phi}(x,\widetilde{y})=\Phi(x,y)$ to obtain
\begin{equation}\label{key.2.2}\begin{split}
\widetilde{\Phi}(x,\widetilde{y})
=\frac{(1-x)(9-\widetilde{y}^2)}{6(x^2\widetilde{y}-3)}
-\frac{1}{18}(x^2 \widetilde{y}+\widetilde{y}-3)(3-3x+\widetilde{y}+x\widetilde{y}), 
\end{split}
\end{equation}
for $(x,y)\in [0,1]\times (0,1]$. For fixed $\widetilde{y}$, it is easy to check $\frac{\partial \widetilde{\Phi}}{\partial x}\geq 0$, then $\widetilde{\Phi}(x,\widetilde{y})\geq \widetilde{\Phi}(0,\widetilde{y})=0$, which completes the proof of the lemma.
\end{proof}
\begin{lemma}\label{lem.2.2}
For weights $\omega_m$ defined in {\rm (\ref{key.1})}, there hold {\rm (i)} $\omega_0>0$, $\omega_m<0, m\geq 1$, {\rm (ii)} $\omega_m>\omega_{m-1}$ $m\geq 2$, and {\rm (iii)} $\displaystyle\sum_{m=0}^{n}\omega_m>0$ for any $n\geq 0$.
\end{lemma}
\begin{proof}
For (i), the correctness for $m=0,1$ can be checked directly.
For $m=2$, (\ref{key.3}) leads to
\begin{equation}\label{key.3.1}\begin{split}
\omega_2=\frac{\alpha}{\alpha+1}\bigg(\frac{1}{\alpha}-\alpha\bigg)\omega_1<0.
\end{split}
\end{equation}
For $m\geq 3$, we first show the following fact
\begin{equation}\label{key.4}
\frac{\omega_{m-1}}{\omega_{m-2}}\geq \sigma_{m-1}>\frac{\frac{1-\alpha}{2\alpha}(m-2)}{\frac{m-1}{\alpha}-\alpha},\quad \text{where~}
\sigma_m=\frac{\frac{1-\alpha}{2\alpha}(m-1)}{\frac{m}{\alpha}-\alpha}\frac{2}{1+\alpha}\big(1+\frac{1}{m}\big),
\end{equation}
then by (\ref{key.3}) we can get
\begin{equation}\label{key.4.1}
\omega_m=\frac{2\alpha\omega_{m-2}}{m(\alpha+1)}\bigg[\frac{1}{\alpha}(m-1)-\alpha\bigg]\bigg[\frac{\omega_{m-1}}{\omega_{m-2}}-\frac{\frac{1-\alpha}{2\alpha}(m-2)}{\frac{1}{\alpha}(m-1)-\alpha}\bigg]<0
\end{equation}
provided that $\omega_{m-2}<0$.
In fact, the inequality (\ref{key.4}) can be derived by induction arguments.
Note first that $\sigma_{m-1}>\frac{\frac{1-\alpha}{2\alpha}(m-2)}{\frac{m-1}{\alpha}-\alpha}$ always holds for $m\geq 3$.
For $m=3$, (\ref{key.4}) is straightforward, as
\begin{equation}\label{key.4.2}
\frac{\omega_2}{\omega_1}=1-\alpha\geq \sigma_2=\frac{3}{2}\frac{1-\alpha}{(2-\alpha^2)(1+\alpha)}.
\end{equation}
Then, for the $m+1$ case, we derive by the recursive formula (\ref{key.3}) that
\begin{equation}\begin{split}\label{key.5}
\frac{\omega_m}{\omega_{m-1}}&\geq \frac{2\alpha}{m(\alpha+1)}\bigg(\frac{m-1}{\alpha}-\alpha-\frac{1-\alpha}{2\alpha}\frac{m-2}{\sigma_{m-1}}\bigg)
\\
&=\frac{2\alpha}{m(\alpha+1)}\bigg(\frac{m-1}{\alpha}-\alpha\bigg)\bigg(1-\frac{1+\alpha}{2}\frac{m-1}{m}\bigg),
\end{split}\end{equation}
where the assumption has been used within the first inequality.
By Lemma \ref{lem.2.1} with $x,y$ replaced with $\alpha$ and $m$, respectively, we have $\frac{\omega_m}{\omega_{m-1}}\geq \sigma_m$ and complete the induction.
\par
For (ii), $\omega_2>\omega_1$ can be directly derived by (\ref{key.3.1}) and the fact $\omega_m<0$ for all $m\geq 1$.
We assume $\frac{\omega_{m-1}}{\omega_{m-2}}<1$ for $m\geq 3$ and again by the recursive formula (\ref{key.3}) we have
\begin{equation}\label{key.6}\begin{split}
\frac{\omega_m}{\omega_{m-1}}&=\frac{2\alpha}{m(\alpha+1)}\bigg[\frac{1}{\alpha}(m-1)-\alpha-\frac{1-\alpha}{2\alpha}(m-2)\frac{\omega_{m-2}}{\omega_{m-1}}\bigg]
\\&\leq
\frac{2\alpha}{m(\alpha+1)}\bigg[\frac{1}{\alpha}(m-1)-\alpha-\frac{1-\alpha}{2\alpha}(m-2)\bigg]
\\&=
1-\frac{2\alpha}{m}<1.
\end{split}\end{equation}
\par
For (iii), obviously we have $\displaystyle\sum_{m=0}^{\infty}\omega_m=0$, then by (i) there holds $\displaystyle\sum_{m=0}^{n}\omega_m>\sum_{m=0}^{\infty}\omega_m$ which completes the proof of the lemma.
\end{proof}
\section{Energy decay and maximum-principle}\label{sec.Energymax}
We first formulate the fully discrete scheme for model (\ref{Int.1}) by adopting finite difference methods in spatial direction and the SFTR-$\frac{1}{2}$ formula in temporal direction for the fractional derivative.
The space $\Omega$ is divided by grid points $a=x_0<x_1<x_2<\dots<x_{M}=b$, $a=y_0<y_1<y_2<\dots<y_{M}=b$, where $x_k=y_k=a+kh$, $h=\frac{b-a}{M}$.
Set $\psi_{j,k}=\psi(x_j,y_k)$.
Take the grid function space $V_h=\{v=\{v_{j,k}\}|v_{j,k}\in \mathbb{R}, 1\leq j,k\leq M\}$ with norms:
\[
\|v\|_\infty=\max_{1\leq j,k\leq M}|v_{j,k}|,\quad
\|v\|^2=h^2\sum_{1\leq j,k\leq M}|v_{j,k}|^2.
\]
For any $u,v \in V_h$, define the discrete $L^2$ inner product as
\[
(u,v)=h^2\sum_{1\leq j,k \leq M}u_{j,k}v_{j,k}.
\]
Introduce the following centered difference formula, for $1\leq j,k\leq M$:
\begin{equation}\label{com.1}\begin{split}
\delta_x^2 \psi_{j,k}=\frac{\psi_{j+1,k}-2\psi_{j,k}+\psi_{j-1,k}}{h^2},&\quad
\delta_y^2 \psi_{j,k}=\frac{\psi_{j,k+1}-2\psi_{j,k}+\psi_{j,k-1}}{h^2}.
\end{split}
\end{equation}
It is a well known result that for sufficiently smooth function $\psi$, there hold
\begin{equation}\label{com.3}\begin{split}
\delta_x^2 \psi_{j,k}&=\partial_x^2\psi_{j,k}+O(h^2),\quad
\delta_y^2 \psi_{j,k}=\partial_y^2\psi_{j,k}+O(h^2).
\end{split}
\end{equation}
Define $\phi^{n-\frac{1}{2}}=\frac{\phi^n+\phi^{n-1}}{2}$ and let $u^n_{j,k}=u(x_j,y_k,t_n)$.
Take $U^n_{j,k}$ as the approximation to $u^n_{j,k}$, then for problem (\ref{Int.1}) at grid point $(x_j,y_k,t_{n-\frac{1}{2}})$, the second-order difference scheme for $1\leq j,k\leq M$, $1\leq n\leq N$ reads
\begin{equation}\label{com.4}\begin{split}
\tau^{-\alpha}\sum_{m=0}^{n}\omega_m (U^{n-m}_{j,k}-U^0_{j,k})
=\varepsilon^2(\delta_x^2+\delta_y^2)U^{n-\frac{1}{2}}_{j,k}
-f^{n-1,n}_{j,k},
\end{split}
\end{equation}
where $U^0_{j,k}=u_0(x_j,y_k)$ and the nonlinear term $f^{n-1,n}_{j,k}$ takes the form
\begin{equation}\label{com.4.1}\begin{split}
f^{n-1,n}_{j,k}=\frac{1}{3}(U^n_{j,k})^{3}+\frac{1}{2}(U^{n-1}_{j,k})^{2}(U^n_{j,k})
+\frac{1}{6}(U^{n-1}_{j,k})^{3}-\frac{1}{2}(U^n_{j,k}+U^{n-1}_{j,k}),
\end{split}
\end{equation}
which immediately leads to the key property (see Appendix A in \cite{LiaoTangZhou1}),
\begin{equation}\label{com.4.2}\begin{split}
(U_{j,k}^n-U_{j,k}^{n-1})f_{j,k}^{n-1,n}\geq F(U_{j,k}^n)-F(U_{j,k}^{n-1}).
\end{split}
\end{equation}
\par
To formulate the related matrix form, for any $v \in V_h$, define $\boldsymbol v \in \mathbb{R}^{M^2}$ by
\[
\boldsymbol v={\rm vec} (v):=[v_{1,1},v_{2,1},\cdots,v_{M,1},v_{1,2},\cdots,v_{M,2},\cdots,v_{1,M},\cdots,v_{M,M}]^T,
\]
with $\|\boldsymbol v\|_\infty=\|v\|_{\infty}$, $\|\boldsymbol v\|=\|v\|$, and
$(\boldsymbol u,\boldsymbol v)=(u,v)$ for any $u,v \in V_h$.
The matrix form of (\ref{com.4}) then reads
\begin{equation}\label{com.5}\begin{split}
\tau^{-\alpha}\sum_{m=0}^{n}\omega_m (\boldsymbol U^{n-m}-\boldsymbol U^0)
=\varepsilon^2 \Delta_{h} \boldsymbol U^{n-\frac{1}{2}}
-\boldsymbol f^{n-1,n},\quad n\geq 1,
\end{split}
\end{equation}
where $\Delta_{h}=\frac{1}{h^2}[I\otimes D+D\otimes I]$ with $\otimes$ standing for the operation of Kronecker tensor product of matrices, and $D$ is a symmetric matrix:
\[
D=
\begin{pmatrix}
  -2 & 1 & 0 & \cdots & 1 \\
  1 & -2 & 1 & \cdots & 0 \\
  \vdots & \ddots & \ddots & \ddots & \vdots \\
  0 & \cdots & 1 & -2 & 1 \\
  1 & \cdots & 0 & 1 & -2
\end{pmatrix}_{M\times M}.
\]
\par
To show the unique solvability of the nonlinear scheme (\ref{com.5}), we reformulate it as
\begin{equation}\label{com.6}\begin{split}
\Upsilon_n \boldsymbol U^n+\frac{1}{3}(\boldsymbol U^n)^{\circ3}=\mathcal{G}(\boldsymbol U^0,\boldsymbol U^1,\cdots,\boldsymbol U^{n-1}),
\end{split}
\end{equation}
where the matrix $\Upsilon_n={\rm Diag}\big(\tau^{-\alpha}\omega_0-\frac{1}{2}+\frac{1}{2}(\boldsymbol U^{n-1})^{\circ2}\big)-\frac{\varepsilon^2}{2}\Delta_h$, and $\mathcal{G}$ stands for the vector
\begin{equation*}\begin{split}
&\mathcal{G}(\boldsymbol U^0,\boldsymbol U^1,\cdots,\boldsymbol U^{n-1})\\
=&
\bigg(\frac{1}{2}I+\frac{\varepsilon^2}{2}\Delta_h\bigg)\boldsymbol U^{n-1}
-\frac{1}{6}(\boldsymbol U^{n-1})^{\circ3}
-\tau^{-\alpha}\sum_{m=1}^{n-1}\omega_m\boldsymbol U^{n-m}
+\tau^{-\alpha}\bigg(\sum_{m=0}^{n-1}\omega_m\bigg)\boldsymbol U^0.
\end{split}
\end{equation*}
Here for the powers with a symbol $\circ$ we mean the calculation is performed elementwisely.
It is a well known result that (\ref{com.6}) admits a unique solution so long as the matrix $\Upsilon_n$ is positive definite \cite{LiaoTangZhou1,LiaoTangZhou,HouTang01}, which is satisfied if we limit $\tau$ such that $\tau^{-\alpha}\omega_0-\frac{1}{2}>0$ as concluded in the following lemma.
\begin{lemma}\label{lem.3.1}
The nonlinear difference scheme (\ref{com.4}) or (\ref{com.5}) is uniquely solvable provided $0<\tau<2^{\frac{1}{\alpha}}\frac{2\alpha}{\alpha+1}$.
\end{lemma}
\subsection{Energy decay}
The discrete energy $E_h^n$ which is an analogue to $E[u](t_n)$ in (\ref{Int.0}) can be defined by
\begin{equation}\label{ene.1}
E_h^n=\frac{\varepsilon^2}{2}\|\nabla_h \boldsymbol U^n\|^2+\frac{1}{4}\|(\boldsymbol U^n)^{\circ2}-1\|^2,
\end{equation}
where $\nabla_h$ is an approximation to $\nabla$ satisfying $(-\Delta_h \boldsymbol U, \boldsymbol V)=(\nabla_h \boldsymbol U, \nabla_h \boldsymbol V)$.
\par
The following sequence is important in deriving the energy decay property.
Define $\{\theta_m\}_{m=0}^\infty$ satisfying
\begin{equation}\label{ene.2}
\sum_{m=0}^{\infty}\theta_m\xi^m=\frac{1-\xi}{\omega(\xi)}=:\theta(\xi),\quad\text{for } |\xi|<1.
\end{equation}
\begin{lemma}\label{lem.3.1.1}
The sequence $\{\theta_m\}_{m=0}^\infty$ defined in (\ref{ene.2}) can be calculated by the recurrence formula:
\begin{equation}\label{ene.3}
\theta_m=\frac{2\alpha}{m(1+\alpha)}\bigg\{\bigg[\frac{1}{\alpha}(m-1)-\frac{1-\alpha}{2\alpha}(2\alpha+1)\bigg]\theta_{m-1}+\frac{1-\alpha}{2\alpha}(3-m)\theta_{m-2}\bigg\}, 
\end{equation}
for $m\geq 2$, with the starting values  $\theta_0=\big(\frac{\alpha+1}{2\alpha}\big)^\alpha$,
$\theta_1=\frac{(\alpha-1)(2\alpha+1)}{\alpha+1}\theta_0$.
\end{lemma}
\begin{proof}
By the definition (\ref{ene.2}), $\theta(\xi)$ has the explicit form
\begin{equation}\label{ene.4}
\theta(\xi)=(1-\xi)^{1-\alpha}\bigg[\frac{1}{2}+\frac{1}{2\alpha}+\bigg(\frac{1}{2}-\frac{1}{2\alpha}\bigg)\xi\bigg]^\alpha.
\end{equation}
Take $p_1(\xi)=\frac{1}{2}+\frac{1}{2\alpha}+\big(\frac{1}{2}-\frac{1}{2\alpha}\big)\xi$ and $p_2(\xi)=1-\xi$, then the formula (\ref{ene.3}) follows from Theorem A.1 of \cite{YinLiuLiZhang333} readily.
\end{proof}
\par
The weights $\theta_m$ actually hold a similar property as that stated in Lemma \ref{lem.2.2} for $\omega_k$, which is determined essentially by the nonnegativity of a special designed function
\begin{equation}\begin{split}\label{ene.5}
\Psi(x,y)&=\bigg(\frac{y-1}{y}-\frac{(1-x)(2x+1)}{2y}\bigg)
\bigg(1-\frac{y-1}{y}\frac{1+x}{2}\bigg)
\\&\quad-\frac{(1-x)(y-2)}{2y-(1-x)(2x+1)}\bigg(1+\frac{1}{y}\bigg).
\end{split}\end{equation}
\begin{lemma}\label{lem.3.1.2}
The function $\Psi(x,y)$ defined in (\ref{ene.5}) is nonnegative for $(x,y)\in [0,1]\times [3,\infty)$.
\end{lemma}
\begin{proof}
Let $y=3/\widetilde{y}$ for $\widetilde{y} \in (0,1]$, and $\widetilde{\Psi}(x,\widetilde{y})=\Psi(x,3/\widetilde{y})$, then one can get
\begin{equation}\label{ene.6}\begin{split}
\widetilde{\Psi}(x,\widetilde{y}) 
&=
\frac{(\widetilde{y}+3)(3-2\widetilde{y})(x-1)}{3\widetilde{y}(2x+1)(x-1)+18}
\\&\quad-\bigg(\frac{(3-\widetilde{y})(x+1)}{6}-1\bigg)
\bigg(\frac{3-\widetilde{y}}{3}+\frac{\widetilde{y}(2x+1)(x-1)}{6}\bigg). \nonumber
\end{split}\end{equation}
Some tedious but simple calculations show that $\frac{\partial \widetilde{\Psi}}{\partial \widetilde{y}}\geq 0$, and thus $\Psi(x,y)=\widetilde{\Psi}(x,\widetilde{y})\geq \widetilde{\Psi}(x,0)=0$.
\end{proof}
\begin{lemma}\label{lem.3.1.3}
For weights $\theta_m$ defined in (\ref{ene.2}), there hold {\rm (i)} $\theta_0>0$, $\theta_m<0, m\geq 1$, and {\rm (ii)} $\displaystyle\sum_{m=0}^{n}\theta_m>0$ for any $n\geq 0$.
\end{lemma}
\begin{proof}
We first prove  {\rm (i)}. The cases for $m=0, 1$ can be verified directly.
For $m=2$, by (\ref{ene.3}) we have $\theta_2=-\frac{2\alpha^3(1-\alpha)}{(1+\alpha)^2}\theta_0<0$.
For $m=3$, similarly, one gets $\theta_3=\frac{2\alpha}{3(1+\alpha)}\big(\alpha+\frac{3-\alpha}{2\alpha}\big)\theta_2<0$.
For $m\geq 4$, we first demonstrate that
\begin{equation}\label{theta.1}
\frac{\theta_{m-1}}{\theta_{m-2}}\geq \delta_{m-1}>\frac{\frac{1-\alpha}{2\alpha}(m-3)}{\frac{m-1}{\alpha}-\frac{1-\alpha}{2\alpha}(2\alpha+1)},
\end{equation}
where 
\begin{equation*}
\delta_m=\frac{\frac{1-\alpha}{2\alpha}(m-2)}{\frac{m}{\alpha}-\frac{1-\alpha}{2\alpha}(2\alpha+1)}\frac{2}{1+\alpha}\bigg(1+\frac{1}{m}\bigg),
\end{equation*}
with which the relation \eqref{ene.3} yields
\begin{equation}\label{theta.1.1}
\theta_m=\frac{2\alpha\theta_{m-2}}{m(1+\alpha)}\bigg\{\bigg[\frac{1}{\alpha}(m-1)-\frac{1-\alpha}{2\alpha}(2\alpha+1)\bigg]\frac{\theta_{m-1}}{\theta_{m-2}}-\frac{1-\alpha}{2\alpha}(m-3)\bigg\}<0
\end{equation}
provided that $\theta_{m-2}<0$.
Actually, \eqref{theta.1} can be proved by induction.
For $m=4$, there holds
\begin{equation}\label{theta.1.2}
\frac{\theta_3}{\theta_2}-\delta_3=\frac{1}{3(1+\alpha)}
\frac{4\alpha^4-4\alpha^3+17\alpha^2+7}{2\alpha^2-\alpha+5}>0.
\end{equation}
We assume that  \eqref{theta.1} is valid, then from \eqref{ene.3} one gets,
\begin{equation}\label{theta.2}
\frac{\theta_{m}}{\theta_{m-1}} = \frac{2\alpha}{1+\alpha}\bigg\{\bigg[\frac{1}{\alpha}\frac{m-1}{m}-\frac{1-\alpha}{2\alpha}(2\alpha+1)\frac{1}{m}\bigg] + \frac{1-\alpha}{2\alpha}\frac{3-m}{m}\frac{\theta_{m-2}}{\theta_{m-1}}\bigg\},
\end{equation}
which, by the assumption, leads to (for $m\geq 4$),
\begin{equation}\label{theta.3}
\begin{split}
\frac{\theta_{m}}{\theta_{m-1}} & \geq \frac{2\alpha}{1+\alpha}\bigg\{\bigg[\frac{1}{\alpha}\frac{m-1}{m}-\frac{1-\alpha}{2\alpha}(2\alpha+1)\frac{1}{m}\bigg] + \frac{1-\alpha}{2\alpha}\frac{3-m}{m}\frac{1}{\delta_{m-1}}\bigg\}\\
& = \frac{2\alpha}{1+\alpha}\bigg[\frac{1}{\alpha}\frac{m-1}{m}-\frac{1-\alpha}{2\alpha}(2\alpha+1)\frac{1}{m}\bigg] \bigg( 1-\frac{m-1}{m}\frac{1+\alpha}{2}\bigg)
\\&\geq
\frac{\frac{1-\alpha}{2\alpha}(m-2)}{\frac{m}{\alpha}-\frac{1-\alpha}{2\alpha}(2\alpha+1)}\frac{2}{1+\alpha}\bigg(1+\frac{1}{m}\bigg)
=\delta_m,
\end{split}
\end{equation}
where the last inequality is due to Lemma \ref{lem.3.1.2} with $x, y$ replaced by $\alpha, m$, respectively.
\par
For {\rm (ii)}, it is obvious that $\displaystyle\sum_{m=0}^{\infty}\theta_m=0$, and by (i) we know that any finite sum with index from $0$ to $n$ is positive.
The proof is completed.
\end{proof}
\begin{remark}\label{rem.3.1.1}
The generating function $\theta(\xi)$ in (\ref{ene.4}) in fact produces a difference formula for the Riemann-Liouville fractional differential operator ${}^R\mathcal{\partial}^{1-\alpha}_t=\partial_t\mathcal{I}^\alpha$ at time $t_{n}$ with second-order accuracy.
To see that, on the one hand $\theta_k$ satisfies $\theta_n=O(n^{-(1-\alpha)-1})$ since $\big[\frac{1}{2}+\frac{1}{2\alpha}+\big(\frac{1}{2}-\frac{1}{2\alpha}\big)\xi\big]^\alpha$
is analytic on the closed unite disc $|\xi|\leq 1$ for any fixed $\alpha \in (0,1)$, and on the other hand,
\begin{equation}\label{ene.7}\begin{split}
\tau^{\alpha-1}\theta(e^{-\tau})&=\bigg(1+\frac{\alpha-1}{2}\tau+\frac{3\alpha^2-7\alpha+4}{24}\tau^2+O(\tau^3)\bigg)
\\&\quad
\times\bigg(1+\frac{1-\alpha}{2}\tau+\frac{(1+\alpha^2)(\alpha-1)}{8\alpha}\tau^2+O(\tau^3)\bigg) \nonumber
\\
&=
1+O(\tau^2). 
\end{split}\end{equation}
Then, the convolution quadrature theory \cite{Lubich1} guarantees $\tau^{\alpha-1}\displaystyle\sum_{m=0}^{n}\theta_{m} \phi^{n-m}$ approximates ${}^R\mathcal{\partial}^{1-\alpha}_t \phi(t_n)$ with second-order accuracy under mild conditions.
\end{remark}
\begin{remark}\label{rem.3.1.3}
We emphasize that the properties revealed in Lemma \ref{lem.2.2} and Lemma \ref{lem.3.1.3} are essential to the following analysis of discrete energy decay and maximum-principle preservation.
The difficulties in the proof for these two lemmas lie in the subtle design of (\ref{key.4}) and (\ref{theta.1}), i.e., the choice of $\sigma_m$ and $\delta_m$, both of which lead to nonnegative functions described in Lemma \ref{lem.2.1} and Lemma \ref{lem.3.1.2}, respectively.
\end{remark}
\par
To discretize the integral part within (\ref{Int.1.3}) based on the weights of SFTR-$\frac{1}{2}$, we introduce another sequence $\{\vartheta_m\}_{m=0}^\infty$ defined as
\begin{equation}\label{ene.8}\begin{split}
\sum_{m=0}^{\infty}\vartheta_m\xi^m=\frac{\theta(\xi)}{1-\xi}=\frac{1}{\omega(\xi)}=:\vartheta(\xi)\quad \text{for } |\xi|<1.
\end{split}\end{equation}
Then, by Cauchy product of series $\displaystyle\sum_{m=0}^{\infty}\theta_m \xi^m$ and $\displaystyle\sum_{m=0}^\infty \xi^m$, i.e.,
\[
\bigg(\displaystyle\sum_{m=0}^{\infty}\theta_m \xi^m\bigg)
\bigg(\displaystyle\sum_{m=0}^\infty \xi^m\bigg)
=\sum_{m=0}^\infty\bigg(\sum_{s=0}^m \theta_s\bigg)\xi^m,
\]
one gets $\vartheta_m=\displaystyle\sum_{s=0}^{m}\theta_s$, followed by the fact
\begin{equation}\label{ene.9}
\vartheta_0>\vartheta_1>\cdots>\vartheta_m>\cdots>0,
\end{equation}
thanks to Lemma \ref{lem.3.1.3}.
By similar arguments stated in Remark \ref{rem.3.1.1}, we have
\begin{equation}\label{ene.10}
\vartheta_n=O(n^{\alpha-1})\quad \text{and }\quad
\tau^\alpha e^{-\frac{1}{2}\tau}\vartheta(e^{-\tau})=1+O(\tau^2).
\end{equation}
Then, the shifted convolution quadrature theory \cite{YinLiuLiZhang444} indicates that, under mild conditions, $\tau^\alpha\displaystyle\sum_{m=0}^{n}\vartheta_m \phi^{n-m}$ approximates $\mathcal{I}_t^\alpha \phi(t_{n+\frac{1}{2}})$ with second-order accuracy.
Now, with the help of the sequence $\{\vartheta_m\}_{m=0}^\infty$, we define the discrete counterpart of $E_c[u]$ in (\ref{Int.1.3}) by
\begin{equation}\label{ene.10.1}
\mathcal{E}^n_h=E_h^n+\frac{\tau^\alpha}{2}\sum_{s=1}^{n}\vartheta_{n-s}\|\boldsymbol V^{s-\frac{1}{2}}\|^2,\quad n\geq 1, \text{ and }\mathcal{E}^0_h=E_h^0,
\end{equation}
where $\boldsymbol V^{n-\frac{1}{2}}:=\varepsilon^2\Delta_h \boldsymbol U^{n-\frac{1}{2}}-\boldsymbol f^{n-1,n}$ is an approximation to $\frac{\delta E}{\delta u}$ at time $t_{n-\frac{1}{2}}$.
\par
Based on Lemma \ref{lem.3.1.3} and the monotonicity of $\vartheta_k$ in (\ref{ene.9}), the following energy decay analysis can be carried out in a standard manner \cite{LiaoTangZhou1}.
We present the details of the analysis within our setting for integrity.
\begin{theorem}\label{thm.3.1.1}
The discrete energy $\mathcal{E}^n_h$ defined in (\ref{ene.10.1}) satisfies the decay property $\mathcal{E}^n_h\leq \mathcal{E}^{n-1}_h$, $n\geq 1$.
\end{theorem}
\begin{proof}
Multiply both hand sides of (\ref{com.5}) by $\tau^{\alpha-1}\theta_{n-s}$ after replacing $n$ with $s$ in (\ref{com.5}), and sum the index $s$ from $1$ to $n$ to obtain
\[
\tau^{-1}\sum_{s=1}^n \theta_{n-s}\sum_{m=0}^s \omega_m(\boldsymbol U^{s-m}-\boldsymbol U^{0})
=\tau^{\alpha-1}\sum_{s=1}^{n}\theta_{n-s}\boldsymbol V^{s-\frac{1}{2}}.
\]
Rearranging the terms in the left hand side of the above equation, one gets
\begin{equation*}\begin{split}
&\tau^{-1}\sum_{s=1}^n \theta_{n-s}\sum_{m=0}^s \omega_m(\boldsymbol U^{s-m}-\boldsymbol U^{0})
=\tau^{-1}\sum_{s=0}^n \theta_{n-s}\sum_{m=0}^s \omega_{s-m}(\boldsymbol U^{m}-\boldsymbol U^{0})
\\
=&\tau^{-1}\sum_{m=0}^n (\boldsymbol U^{m}-\boldsymbol U^{0})  \sum_{s=m}^n \omega_{s-m} \theta_{n-s}
=\tau^{-1}\sum_{m=0}^n (\boldsymbol U^{m}-\boldsymbol U^{0})  \sum_{s=0}^{n-m} \omega_{s} \theta_{n-s-m}
\\
=&\tau^{-1}\sum_{m=0}^n (\boldsymbol U^{n-m}-\boldsymbol U^{0})  c_m
\end{split}\end{equation*}
where $c_m=\sum_{s=0}^{m} \omega_{s} \theta_{m-s}$ is the $m$-th coefficient of the series of the product $\omega(\xi)\theta(\xi)$, which by (\ref{ene.2}) yields $c_0=1, c_1=-1$ and $c_k=0$ for $k>1$.
We then have
\begin{equation}\label{ene.12}
\frac{\boldsymbol U^n-\boldsymbol U^{n-1}}{\tau}=\tau^{\alpha-1}\sum_{s=1}^{n}\theta_{n-s}\boldsymbol V^{s-\frac{1}{2}}.
\end{equation}
Taking the $L^2$ inner product of (\ref{ene.12}) with respect to $\boldsymbol V^{n-\frac{1}{2}}$, we get
\begin{equation}\label{ene.13}
\bigg(\boldsymbol V^{n-\frac{1}{2}},\sum_{s=1}^{n}\theta_{n-s}\boldsymbol V^{s-\frac{1}{2}}\bigg)
=\tau^{-\alpha}\big(\boldsymbol U^n-\boldsymbol U^{n-1},\varepsilon^2\Delta_h \boldsymbol U^{n-\frac{1}{2}}-\boldsymbol f^{n-1,n}\big).
\end{equation}
Assume $\boldsymbol V^s={\rm vec}(v^s)$ for some $v^s=\{v_{j,k}^s\}\in V_h$.
Since $\theta_s=\vartheta_s-\vartheta_{s-1}<0, s\geq 1$ and $\theta_0=\vartheta_0>0$ by Lemma \ref{lem.3.1.3}, there holds
\begin{equation}\label{ene.14}\begin{split}
&\quad v^{n-\frac{1}{2}}_{j,k}\sum_{s=1}^{n}\theta_{n-s}v^{s-\frac{1}{2}}_{j,k}
=\vartheta_0(v^{n-\frac{1}{2}}_{j,k})^2
+\sum_{s=1}^{n-1}\theta_{n-s}\big(v^{s-\frac{1}{2}}_{j,k}\big)\big(v^{n-\frac{1}{2}}_{j,k}\big)
\\&\geq
\vartheta_0(v^{n-\frac{1}{2}}_{j,k})^2
+\frac{1}{2}\sum_{s=1}^{n-1}\theta_{n-s}\big[\big(v^{s-\frac{1}{2}}_{j,k}\big)^2+\big(v^{n-\frac{1}{2}}_{j,k}\big)^2\big]
\\&=
\frac{1}{2}\vartheta_{n-1}\big(v^{n-\frac{1}{2}}_{j,k}\big)^2
+\frac{1}{2}\sum_{s=1}^{n}\vartheta_{n-s}\big(v^{s-\frac{1}{2}}_{j,k}\big)^2
-\frac{1}{2}\sum_{s=1}^{n-1}\vartheta_{n-s-1}\big(v^{s-\frac{1}{2}}_{j,k}\big)^2.
\end{split}\end{equation}
Then, for the left hand side of (\ref{ene.13}) we have
\begin{equation}\label{ene.15}\begin{split}
& \bigg(\boldsymbol V^{n-\frac{1}{2}},\sum_{s=1}^{n}\theta_{n-s}\boldsymbol V^{s-\frac{1}{2}}\bigg)\\
\geq & \frac{1}{2}\vartheta_{n-1}\|\boldsymbol V^{n-\frac{1}{2}}\|^2
+\frac{1}{2}\sum_{s=1}^{n}\vartheta_{n-s}\|\boldsymbol V^{s-\frac{1}{2}}\|^2
-\frac{1}{2}\sum_{s=1}^{n-1}\vartheta_{n-s-1}\|\boldsymbol V^{s-\frac{1}{2}}\|^2.
\end{split}\end{equation}
Moreover, since $\big(\boldsymbol U^n-\boldsymbol U^{n-1},\Delta_h \boldsymbol U^{n-\frac{1}{2}}\big)=-\big(\|\nabla_h \boldsymbol U^n\|^2-\|\nabla_h \boldsymbol U^{n-1}\|^2\big)$, combined with (\ref{com.4.2}), the right hand side of (\ref{ene.13}) can be estimated by
\begin{equation}\label{ene.16}\begin{split}
&\quad\tau^{-\alpha}\big(\boldsymbol U^n-\boldsymbol U^{n-1},\varepsilon^2\Delta_h \boldsymbol U^{n-\frac{1}{2}}-\boldsymbol f^{n-1,n}\big)
\\&\leq
\frac{-\varepsilon^2}{2\tau^\alpha}\big(\|\nabla_h \boldsymbol U^n\|^2-\|\nabla_h \boldsymbol U^{n-1}\|^2\big)
-\frac{h^2}{\tau^{-\alpha}}\sum_{1\leq j,k \leq M}\big(F(U_{j,k}^n)-F(U_{j,k}^{n-1})\big)
\\&=
-\frac{\varepsilon^2}{2\tau^\alpha}\big(\|\nabla_h \boldsymbol U^n\|^2-\|\nabla_h \boldsymbol U^{n-1}\|^2\big)
\\
&\quad-\frac{1}{4\tau^\alpha}\big[\|\big(\boldsymbol U^n\big)^{\circ2}-1\|^2-\|\big(\boldsymbol U^{n-1}\big)^{\circ2}-1\|^2\big],
\end{split}\end{equation}
which combined with (\ref{ene.15}) leads to
\begin{equation}\label{ene.17}\begin{split}
\mathcal{E}^n_h-\mathcal{E}^{n-1}_h+\frac{\tau^\alpha}{2}\vartheta_{n-1}\|\boldsymbol V^{n-\frac{1}{2}}\|^2\leq 0.
\end{split}\end{equation}
This completes the proof of the theorem.
\end{proof}
\begin{remark}\label{rem.3.1.2}
It is notable that when $\alpha \to 1$, $\theta(\xi)=1$ by (\ref{ene.4}) and thus $\vartheta_m=1$ for $m\geq 0$.
In this case, $\mathcal{E}^n_h$ defined in (\ref{ene.10.1}) becomes to
\begin{equation}\label{ene.18}
\mathcal{E}^n_h=E_h^n+\frac{\tau}{2}\sum_{s=1}^{n}\|\boldsymbol V^{s-\frac{1}{2}}\|^2,\quad n\geq 1, \text{ and }\mathcal{E}^0_h=E_h^0,
\end{equation}
and (\ref{ene.17}) degrades into
\begin{equation}\label{ene.19}\begin{split}
\frac{1}{\tau}(E^n_h-E^{n-1}_h)+\|\boldsymbol V^{n-\frac{1}{2}}\|^2\leq 0,
\end{split}\end{equation}
which is compatible with the classical case (\ref{Int.1.1}).
\end{remark}
\subsection{Maximum-principle}\label{sec.von}
The properties of $\omega_m$ proved in Lemma \ref{lem.2.2} and the form of nonlinear term approximation (\ref{com.4.1}) followed by some standard analysis essentially lead to the discrete maximum-principle (see \cite[Theorem 3.2]{LiaoTangZhou1}).
We shall present in the next theorem the outline of the proof to find out the restriction on time-step size $\tau$ for our scheme.
%
\begin{theorem}\label{thm.3.2.1}
The numerical scheme (\ref{com.4}) or (\ref{com.5}) preserves the discrete maximum-principle, i.e., for
$\tau$ satisfying $0<\tau<\min\big\{2^{\frac{1}{\alpha}}\frac{2\alpha}{\alpha+1},(\frac{\alpha h^2}{2\varepsilon^2})^{\frac{1}{\alpha}}(\frac{2\alpha}{\alpha+1})^{\frac{\alpha+1}{\alpha}}\big\}$,
there holds $\|\boldsymbol U^n\|_{\infty}\leq 1$ for all $n \geq 1$
provided $u_0(\boldsymbol x)\in C(\overline{\Omega})$ and $\max_{\boldsymbol x \in \overline{\Omega}}|u_0(\boldsymbol x)|\leq 1$.
\end{theorem}
\begin{proof}
With $\Upsilon_n$ defined in (\ref{com.6}), we rewrite (\ref{com.5}) into
\begin{equation}\begin{split}\label{ene.20}
&\quad\Upsilon_n \boldsymbol U^n+\frac{1}{3}(\boldsymbol U^n)^{\circ3}
\\&=\mathcal{G}_1(\boldsymbol U^{n-1})
+\mathcal{G}_2(\boldsymbol U^{0},\boldsymbol U^{1},\cdots,\boldsymbol U^{n-2})
+\frac{1}{6}\big[3\boldsymbol U^{n-1}-(\boldsymbol U^{n-1})^{\circ3}\big],
\end{split}\end{equation}
where
\begin{equation}\label{ene.21}\begin{split}
\mathcal{G}_1(\boldsymbol U^{n-1})&=\bigg(\frac{\varepsilon^2}{2}\Delta_h-\tau^{-\alpha}\omega_1I\bigg)\boldsymbol U^{n-1},
\\
\mathcal{G}_2(\boldsymbol U^{0},\boldsymbol U^{1},\cdots,\boldsymbol U^{n-2})&=
-\tau^{-\alpha}\sum_{m=2}^{n-1}\omega_m\boldsymbol U^{n-m}
+\tau^{-\alpha}\bigg(\sum_{m=0}^{n-1}\omega_m\bigg)\boldsymbol U^0.
\end{split}\end{equation}
Assume $\|\boldsymbol U^m\|_\infty\leq 1$ for $m=0,1,\cdots,n-1$, and let $G_h=(g_{jk})=\big(\frac{\varepsilon^2}{2}\Delta_h-\tau^{-\alpha}\omega_1I\big)$.
Then, there hold $g_{jk}>0$ for $j\neq k$ and $g_{jj}=-\frac{2\varepsilon^2}{h^2}-\tau^{-\alpha}\omega_1$, and $\sum_{k}g_{jk}=-\tau^{-\alpha}\omega_1$ for any $j$.
If $\tau\leq \big(\frac{\alpha h^2}{2\varepsilon^2}\big)^{\frac{1}{\alpha}}\big(\frac{2\alpha}{\alpha+1}\big)^{\frac{\alpha+1}{\alpha}}$, i.e., $g_{jj}\geq 0$, then
\begin{equation}\label{ene.22}\begin{split}
\|\mathcal{G}_1(\boldsymbol U^{n-1})\|_\infty
\leq\|G_h\|_\infty\|\boldsymbol U^{n-1}\|_\infty
\leq -\frac{1}{2}\tau^{-\alpha}\omega_1(1+\|\boldsymbol U^{n-1}\|_\infty).
\end{split}\end{equation}
Thanks to Lemma \ref{lem.2.2}, the $\mathcal{G}_2$ term can be bounded by $\|\mathcal{G}_2\|_\infty\leq \tau^{-\alpha}(\omega_0+\omega_1)$.
Since $|3z-z^3|\leq 2$ for any $z\in [-1,1]$, then one gets $\|\frac{1}{6}\big[3\boldsymbol U^{n-1}-(\boldsymbol U^{n-1})^{\circ3}\big]\|_\infty\leq \frac{1}{3}$.
Moreover, the unique solvability condition on $\tau$ in Lemma \ref{lem.3.1}, i.e., $\tau<2^{\frac{1}{\alpha}}\frac{2\alpha}{\alpha+1}$ can lead to the following estimates (see \cite[Lemma 3.5]{LiaoTangZhou1}),
\begin{equation}
\begin{split}\label{ene.23}
&\quad\bigg\|\Upsilon_n \boldsymbol U^n+\frac{1}{3}(\boldsymbol U^n)^{\circ3}\bigg\|_\infty
\\&\geq
(\tau^{-\alpha}\omega_0-\frac{1}{2})\|\boldsymbol U^n\|_\infty
+\frac{1}{2}\|\boldsymbol U^{n-1}\|_\infty^2 \|\boldsymbol U^{n}\|_\infty
+\frac{1}{3}\|\boldsymbol U^n\|_\infty^3.
\end{split}
\end{equation}
The above term by term estimates of (\ref{ene.20}) combined with a subtle contradiction argument finally give the desired conclusion.
More details can be found in \cite{LiaoTangZhou1}.
\end{proof}

\begin{remark}
In this paper, our focus is on the structure-preserving properties of the scheme. As for its convergence, here we give some remarks. If the exact solution is sufficiently smooth, it is easy to prove the second-order convergence of the scheme by using the standard energy method. Here we omit the proof because it is more reasonable to think that the solution has a weak singularity at $t=0$. In the next section, we will give a numerical investigation, and the results show that the method has  second order convergence at fixed positive time when solving the tFAC equation \eqref{Int.1}. 
Despite the rigorous arguments of sharp error estimate is still open to us, we would like to give some insights into why SFTR-$\frac{1}{2}$ automatically preserve the optimal accuracy at fixed time in numerical tests.
Actually, for the linear counterpart of (\ref{Int.1}), i.e.,  the linear subdiffusion problem with $f(u)=0$,  rigorous arguments have been carried out in \cite{YinLiuLiZhang222} for SFTR-$\theta$ with $\theta \in (0,\frac{1}{2})$, where the correction term $$(1/2-\theta)\tilde{\Delta}_h u_{0h}$$ is added to the initial time step to compensate for the weak regularity of the solution and to restore the second-order accuracy.
$\tilde{\Delta}_h$ denotes some discrete approximation to the Laplacian.
The correction term can be arbitrarily small so long as the given $\theta$ approaches $\frac{1}{2}$ adequately while the second-order accuracy is still guaranteed, and it is then quite natural to ask what may happen if $\theta$ takes $\frac{1}{2}$, since in this case the correction term vanishes formally.
However, the rigorous theoretical analysis for the limit case $\theta=\frac{1}{2}$ is still an open problem for that the arguments developed for $\theta \in (0,\frac{1}{2})$ can not be simply carried out by setting $\theta=\frac{1}{2}$ directly.
As shown in \cite{YinLiuLiZhang222}, the arguments for $\theta \in (0,\frac{1}{2})$ rely essentially on the solution representation theory where a contour integral is involved.
To be specific, let $U_h^n$ be the approximation to $u(t_n)$ within some space discretization framework such as the finite element  method, then $U_h^n$ can be represented by
\begin{equation}\label{ene.24}
U_h^n=u_{h0}-\frac{1}{2\pi \rm{i}}\int_{\mathcal{C}^\tau}e^{zt_n}\big[
\mu(e^{-z\tau})K(\beta_\tau(e^{-z\tau}))\tilde{\Delta}_h u_{h0}
\big],
\end{equation}
where $K(z)=-z^{-1}(z^\alpha-\tilde{\Delta}_h)^{-1}$, $\beta_\tau(\zeta)=\frac{1-\zeta}{\tau(1+\mu_1\zeta)}\big(\frac{\mu_0}{1-\theta+\theta\zeta}\big)^{\frac{1}{\alpha}}$ and $\mu(\zeta)=\frac{\mu_0^{\frac{1}{\alpha}}(3/2-\theta)\zeta(1-\mu_2\zeta)}{(1+\mu_1\zeta)(1-\theta+\theta\zeta)^{1+\frac{1}{\alpha}}}$ with the coefficients $\mu_0, \mu_1$ and $\mu_2$ defined by $\mu_0=\big(\frac{2\alpha}{\alpha+2\theta}\big)^\alpha$, $\mu_1=\frac{\alpha-2\theta}{\alpha+2\theta}$ and $\mu_2=\frac{1-2\theta}{3-2\theta}$.
\par
To obtain (\ref{ene.24}), it is required (Theorem 4.3 in  \cite{YinLiuLiZhang222}) that both $\beta_\tau(\zeta)$ and $\mu(\zeta)$ are analytic at $\zeta \in \{\zeta \in \mathbb{C}:|\zeta|\leq 1, \zeta \neq 1\}$.
Unfortunately, when $\theta=\frac{1}{2}$, both $\beta_\tau(\zeta)$ and $\mu(\zeta)$ are singular at point $\zeta=-1$, which implies the underlying arguments for $\theta=\frac{1}{2}$ is non-standard and that replacing $\theta$ with $\frac{1}{2}$ directly is not acceptable.
For the nonlinear subdiffusion problem such as the tFAC equation, more technique should be involved.
Despite the difficulties and unknowns in the error estimates, the following numerical simulations are still valuable which, for the one hand suggests SFTR-$\frac{1}{2}$ is superior to the averaged fractional BDF2 on uniform meshes and is comparable to L2-1$_\sigma$ method on graded meshes, and for the other hand convince us by the idea that difference methods proposed at shifted non-integer grids can resolve the initial singularity of fractional problems to some extent.
\end{remark}

\section{Numerical experiments}\label{sec.tests}
We implement several examples in this section to verify the energy decay and maximum-principle preservation, and demonstrate the high efficiency of SFTR-$\frac{1}{2}$ under low solution regularity.
Since the fully discrete scheme (\ref{com.4}) or (\ref{com.5}) is nonlinear, we adopt a simple iteration method with the stopping criterion that $\|U^{(k+1)}-U^{(k)}\|_\infty\leq 10^{-6}$.
\subsection{Verification of energy decay and maximum-principle}
\par
\begin{figure}[htbp]
\centering
\subfigure[]{
\begin{minipage}[t]{0.3\linewidth}
\centering
\includegraphics[width=1\textwidth]{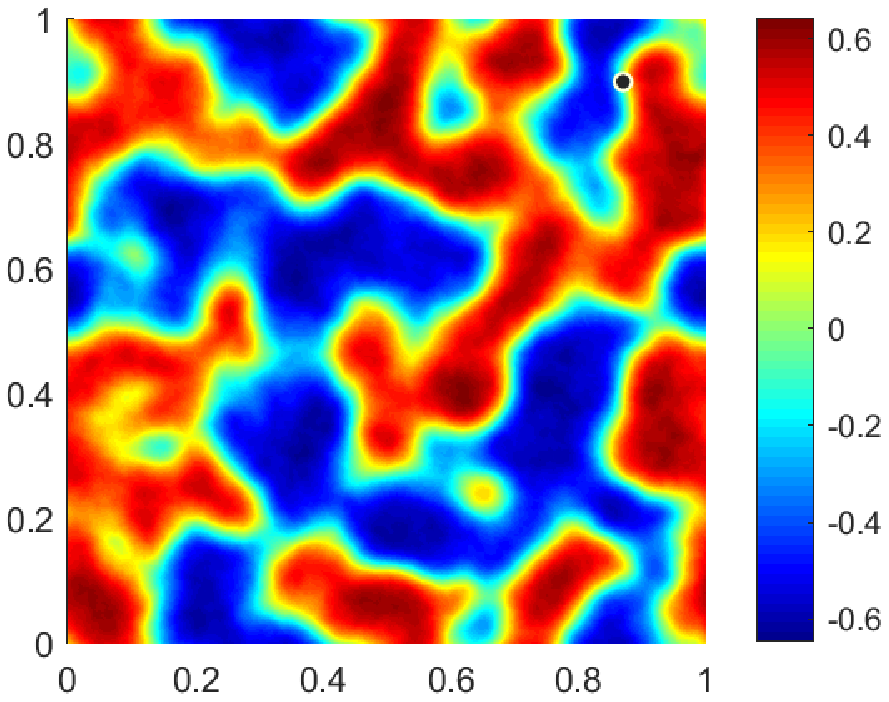}
\end{minipage}%
}%
\subfigure[]{
\begin{minipage}[t]{0.3\linewidth}
\centering
\includegraphics[width=1\textwidth]{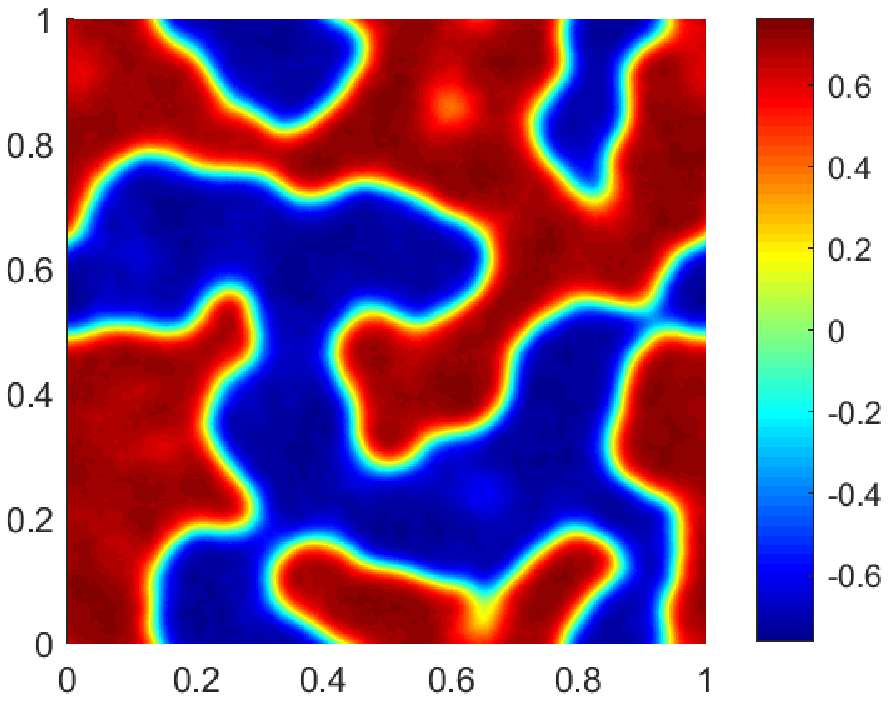}
\end{minipage}%
}%
\subfigure[]{
\begin{minipage}[t]{0.3\linewidth}
\centering
\includegraphics[width=1\textwidth]{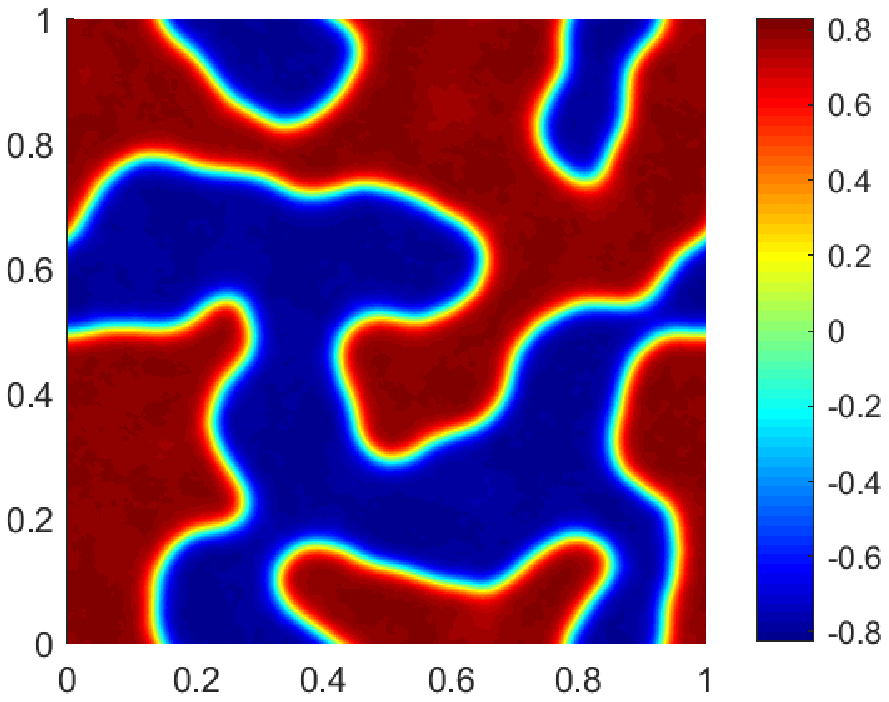}
\end{minipage}%
}%
\\
\subfigure[]{
\begin{minipage}[t]{0.3\linewidth}
\centering
\includegraphics[width=1\textwidth]{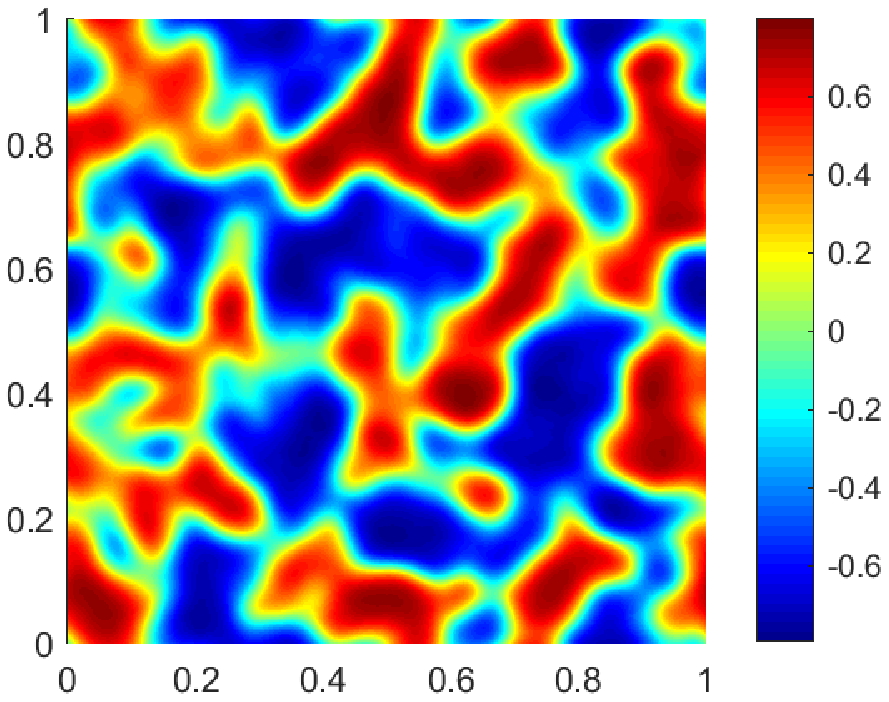}
\end{minipage}%
}%
\subfigure[]{
\begin{minipage}[t]{0.3\linewidth}
\centering
\includegraphics[width=1\textwidth]{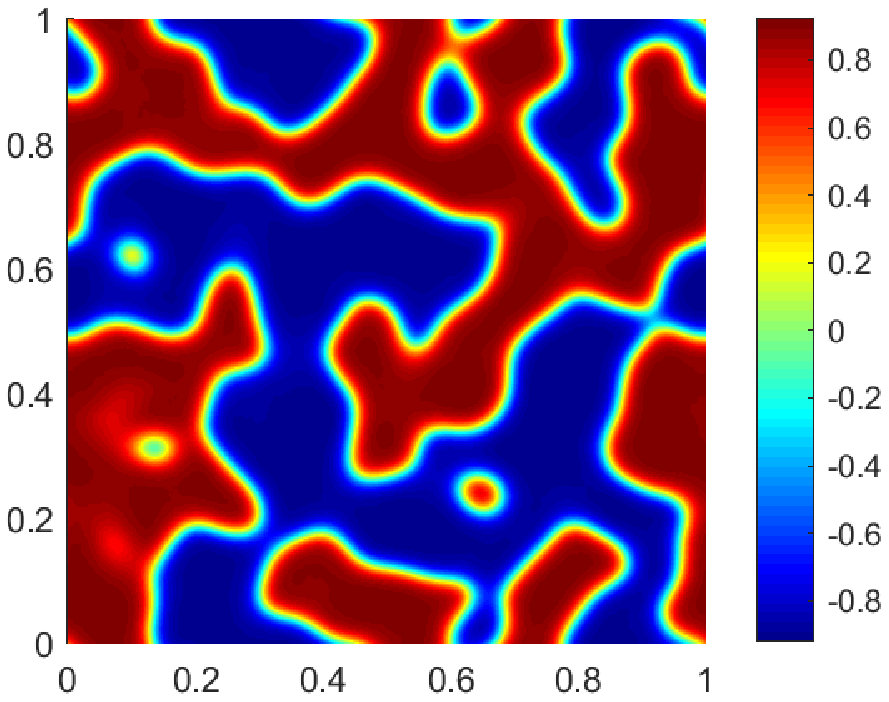}
\end{minipage}%
}%
\subfigure[]{
\begin{minipage}[t]{0.3\linewidth}
\centering
\includegraphics[width=1\textwidth]{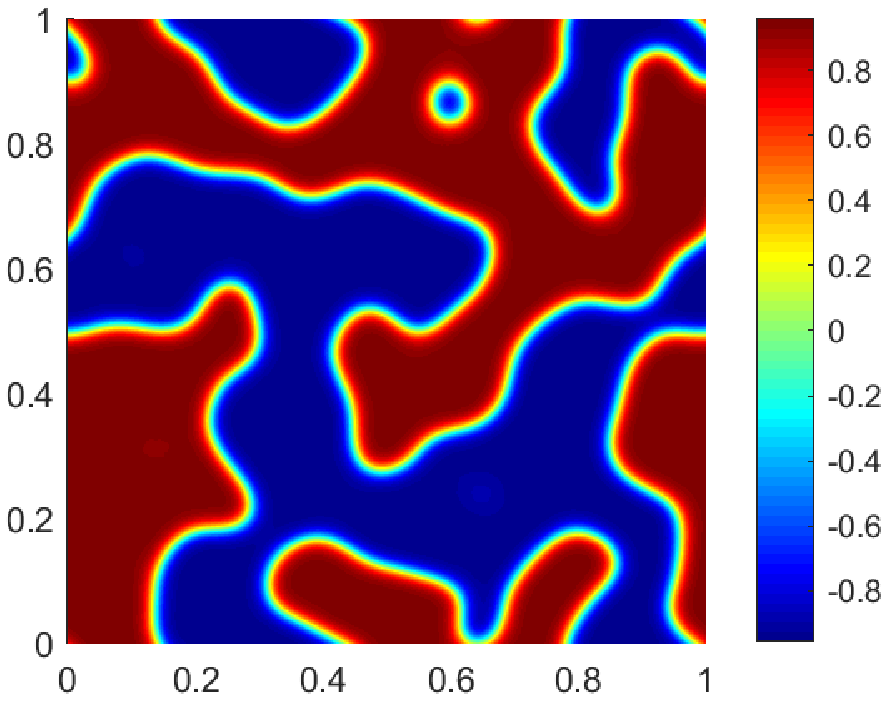}
\end{minipage}%
}%
\\
\subfigure[]{
\begin{minipage}[t]{0.3\linewidth}
\centering
\includegraphics[width=1\textwidth]{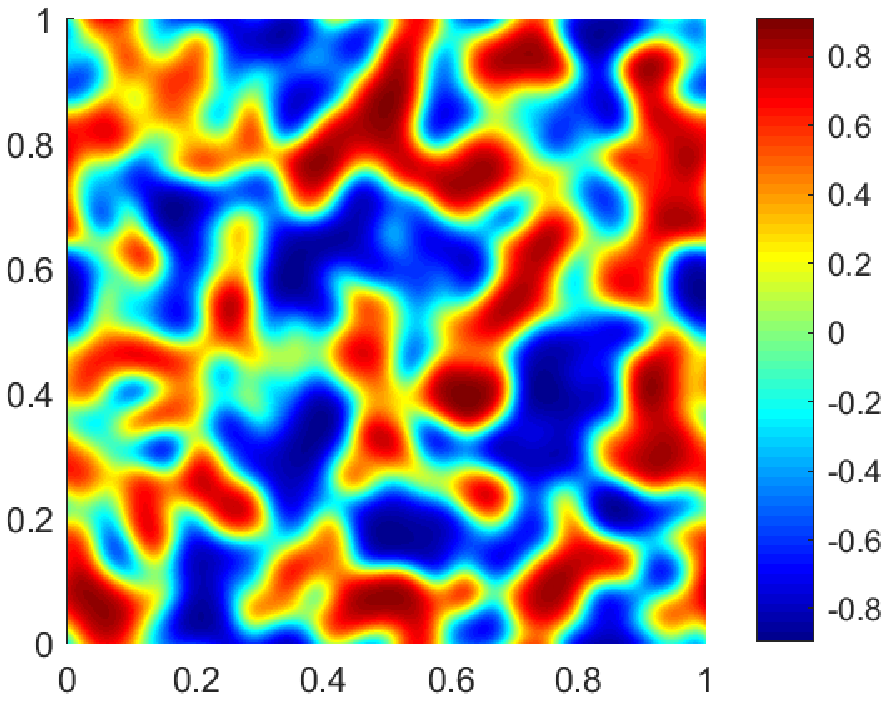}
\end{minipage}%
}%
\subfigure[]{
\begin{minipage}[t]{0.3\linewidth}
\centering
\includegraphics[width=1\textwidth]{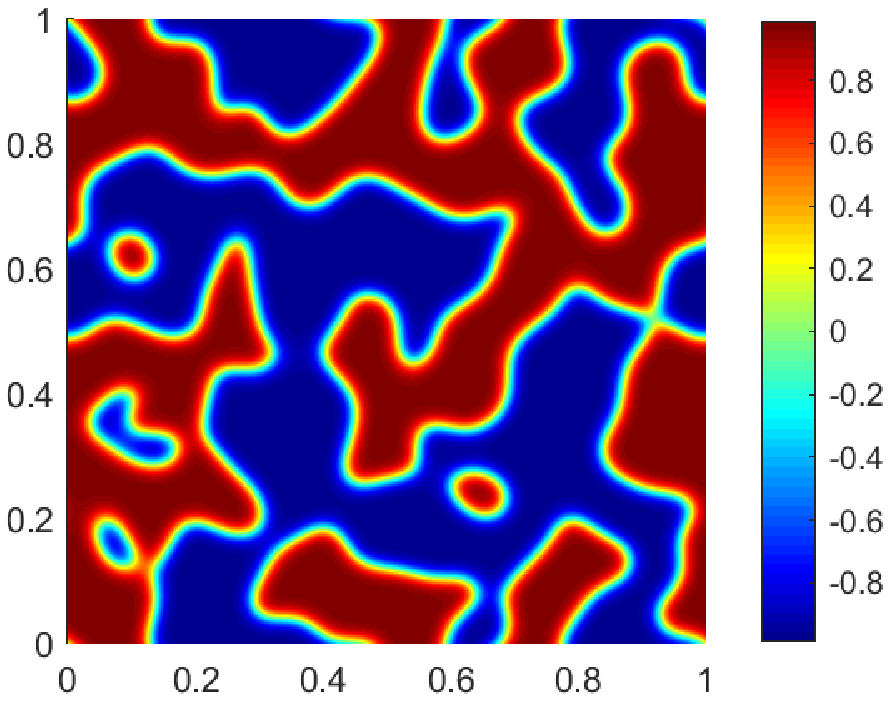}
\end{minipage}%
}%
\subfigure[]{
\begin{minipage}[t]{0.3\linewidth}
\centering
\includegraphics[width=1\textwidth]{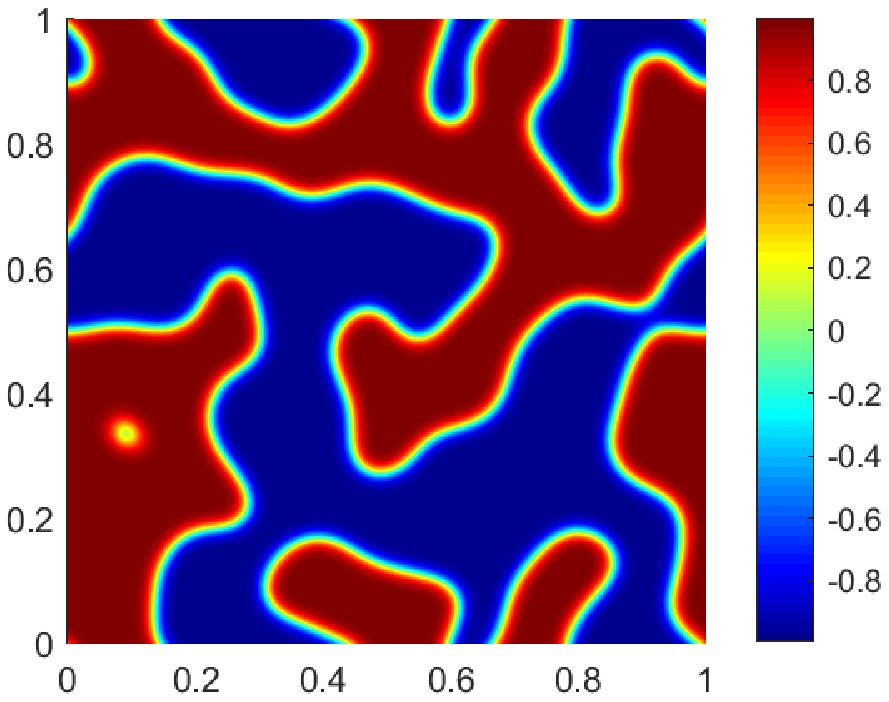}
\end{minipage}%
}%
\centering
\caption{Snapshots of \textit{Example 1} at time $t=5,10,20$ (from left to right) for different $\alpha=0.3,0.6 $ and $ 0.9$ (from top to bottom), respectively.}\label{C1}
\end{figure}
\begin{figure}[htbp]
\centering
\subfigure[]{
\begin{minipage}[t]{0.45\linewidth}
\centering
\includegraphics[width=1\textwidth]{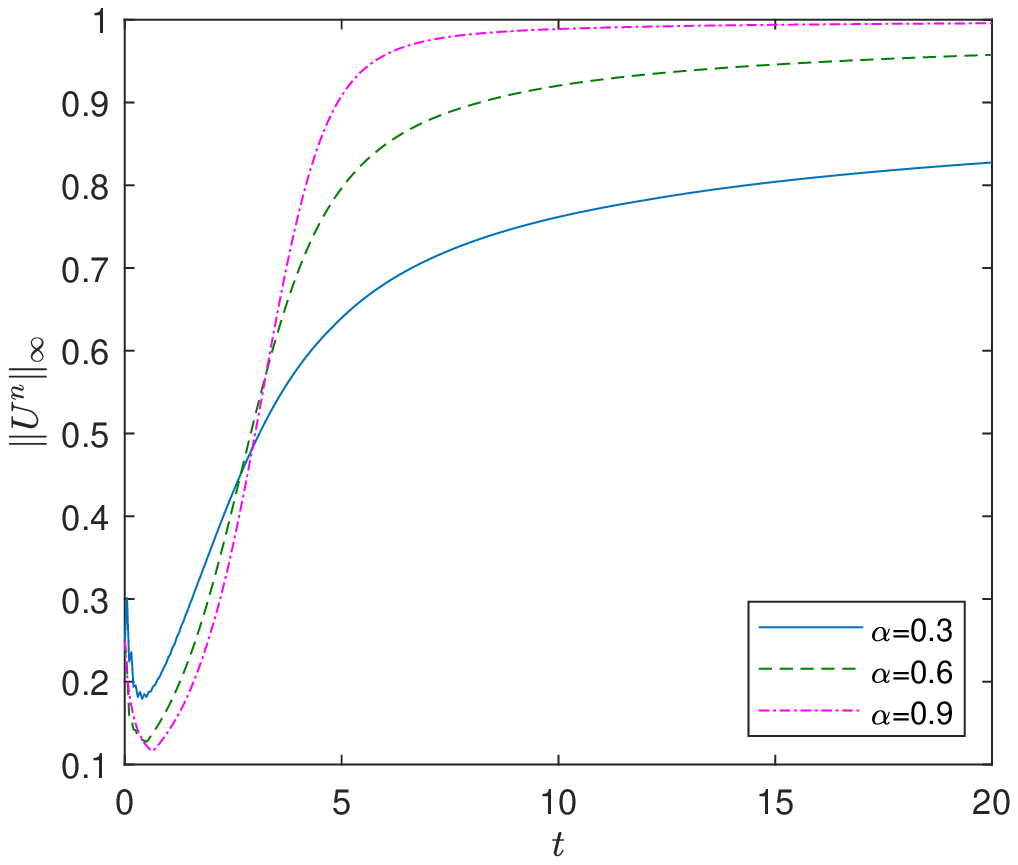}
\end{minipage}%
}%
\subfigure[]{
\begin{minipage}[t]{0.45\linewidth}
\centering
\includegraphics[width=1\textwidth]{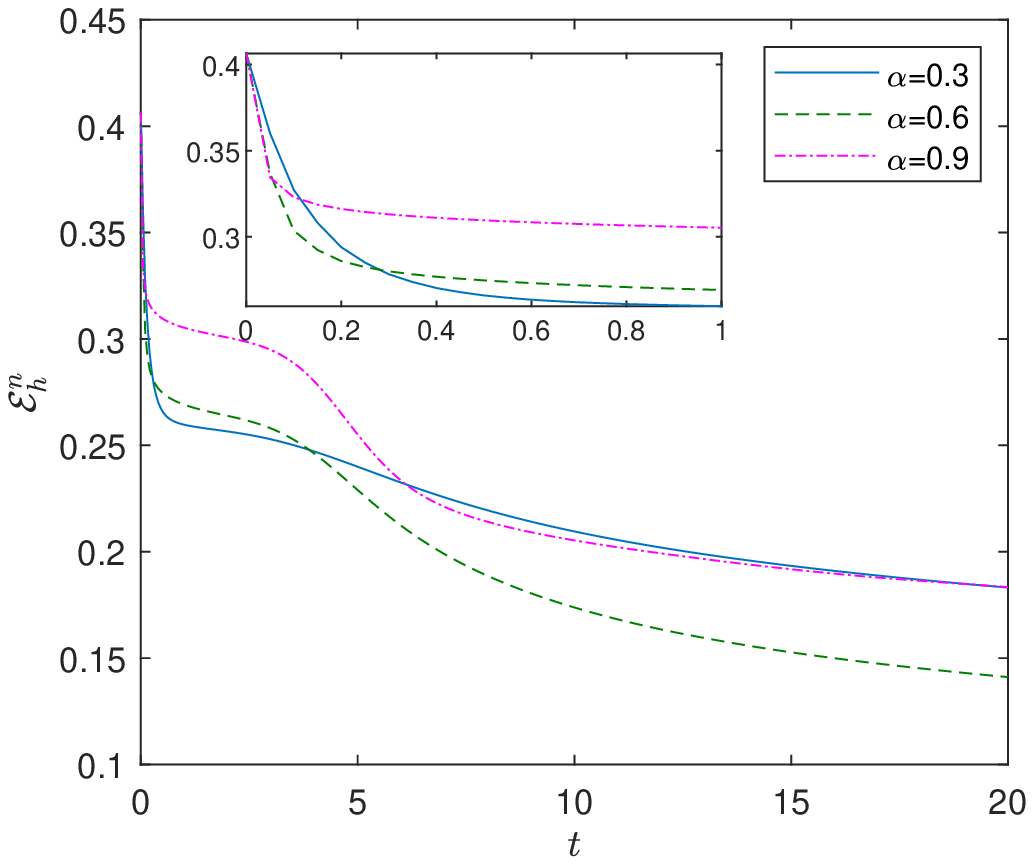}
\end{minipage}%
}%
\centering
\caption{Maximum norm $\|U^n\|_\infty$ (a) and discrete compatible energy $\mathcal{E}_h^n$ (b) in \textit{Example 1} for different $\alpha=0.3,0.6$ and $0.9$, respectively.}\label{C2}
\end{figure}

\textit{Example 1.}\quad
Let $\Omega=(0,1)\times(0,1)$ and $\varepsilon=0.01$ and take the initial condition $u_0(\boldsymbol x)=0.5\times{\rm rand}(\boldsymbol x)-0.25$ where the function ${\rm rand}(\boldsymbol x)$ returns uniformly distributed random numbers in $(0,1)$ and thus $|u_0(\boldsymbol x)|\leq 1$.
The time space mesh is chosen with $\tau=0.05$ and $h=\frac{1}{200}$.
We simulate the coarsening dynamics under the above setting and illustrate some snapshots in Fig.\ref{C1} up to $T=20$ for different fractional differential order $\alpha=0.3,0.6$ and $0.9$, respectively.
The maximum of $|U^n|$ and discrete compatible energy $\mathcal{E}_h^n$ for each time level are depicted in Fig.\ref{C2} (a) and (b), for $\alpha=0.3,0.6,0.9$ respectively.
Obviously, under this random initial condition, the discrete maximum-principle and energy decay property are preserved.

\subsection{High accuracy under low regularity}
We demonstrate in the following two examples the serious impact of low regularity of solution on temporal convergence accuracy by first simulating a model with smooth enough solution (in which case the source term is nonzero), and then exploring the original model and comparing the results between our method and the well known fractional BDF2 (F-BDF2) \cite{Lubich1}.
For both examples, let $\Omega=(0,1)\times(0,1)$, $T=1$ and $\varepsilon=0.005$.
\\\\
\textit{Example 2.}\quad
To approximate $\partial_t^\alpha u$ at time $t_{n-\frac{1}{2}}$ by the F-BDF2, we take the following process,
\begin{equation}\label{tes.1}
\partial_t^\alpha u(t_{n-\frac{1}{2}})\approx
\frac{1}{2}\big[\partial_t^\alpha u(t_{n})+\partial_t^\alpha u(t_{n-1})\big],
\end{equation}
where $\partial_t^\alpha u(t_{n})\approx \tau^{-\alpha}\sum_{m=0}^{n}\widetilde{\omega}_{n-m}(u(t_m)-u(t_0))$ with $\widetilde{\omega}_m$ generated by the function
\[
\widetilde{\omega}(\xi)=\sum_{m=0}^{\infty}\widetilde{\omega}_m\xi^m=\bigg(\frac{3}{2}-2\xi+\frac{1}{2}\xi^2\bigg)^\alpha.
\]
\par
Assume the exact solution is $$u(\boldsymbol x,t)=u(x_1,x_2,t)=\frac{1}{4}(1+t^3)\sin(2\pi x_1)\sin(2\pi x_2).$$
The source term can be derived easily.
Note that such solution is smooth in time from which nice convergence results can be expected.
The temporal rates at time level $t_n$ is obtained by the formula:
\[
\log_2\frac{\text{Error}(\tau;h,t_n)}{\text{Error}(\frac{\tau}{2};h,t_n)},\quad\text{where~~ }
\text{Error}(\tau;h,t_n)=\|U^n-u(\boldsymbol x,t_n)\|.
\]
\par
In Table \ref{tab1}, we report for different $\alpha=0.3,0.6,0.9$ the $L^2$ error at time $t=\frac{T}{2}$ and convergence rates for both SFTR-$\frac{1}{2}$ and F-BDF2, by taking different time-step size $\tau=\frac{1}{20},\frac{1}{40},\frac{1}{80}$ and $\frac{1}{160}$, respectively.
For such smooth solution, the optimal second-order convergence accuracy is obtained which is in line with our expectation.
\begin{table}[]
\centering
\caption{Results for smooth solutions in \textit{Example 2}.}\label{tab1}
{\footnotesize
{\renewcommand{\arraystretch}{1.4}
\begin{tabular*}{\hsize}{@{}@{\extracolsep{\fill}}cllllll@{}}
\toprule
\multirow{2}{*}{$\alpha$} & \multicolumn{1}{c}{\multirow{2}{*}{$\tau$}} & \multicolumn{2}{c}{SFTR-$\frac{1}{2}$}                &  & \multicolumn{2}{c}{F-BDF2}                            \\ \cline{3-4} \cline{6-7}
                          & \multicolumn{1}{c}{}                        & \multicolumn{1}{c}{Errer} & \multicolumn{1}{c}{Rates} &  & \multicolumn{1}{c}{Errer} & \multicolumn{1}{c}{Rates} \\ \hline
\multirow{4}{*}{0.3}      & 1/20                                        & 5.9056E-04                &                           &  & 8.6484E-05                &                           \\
                          & 1/40                                        & 1.6353E-04                & 1.85                      &  & 2.3802E-05                & 1.86                      \\
                          & 1/80                                        & 4.3146E-05                & 1.92                      &  & 6.0214E-06                & 1.98                      \\
                          & 1/160                                       & 1.1200E-05                & 1.95                      &  & 1.2785E-06                & 2.24                      \\ \hline
\multirow{4}{*}{0.6}      & 1/20                                        & 1.5359E-04                &                           &  & 1.5720E-04                &                           \\
                          & 1/40                                        & 4.0232E-05                & 1.93                      &  & 4.2323E-05                & 1.89                      \\
                          & 1/80                                        & 1.0320E-05                & 1.96                      &  & 1.0950E-05                & 1.95                      \\
                          & 1/160                                       & 2.6647E-06                & 1.95                      &  & 2.6695E-06                & 2.04                      \\ \hline
\multirow{4}{*}{0.9}      & 1/20                                        & 4.0839E-05                &                           &  & 2.2212E-04                &                           \\
                          & 1/40                                        & 1.0424E-05                & 1.97                      &  & 5.9128E-05                & 1.91                      \\
                          & 1/80                                        & 2.6237E-06                & 1.99                      &  & 1.5253E-05                & 1.95                      \\
                          & 1/160                                       & 6.7011E-07                & 1.97                      &  & 3.8653E-06                & 1.98                      \\
\bottomrule
\end{tabular*}}}
\end{table}
\\\\
\textit{Example 3.} This example is devoted to the model (\ref{Int.1}), i.e., with zero source term.
Take the initial condition as $u_0(\boldsymbol x)=u_0(x_1,x_2)=\sin(2\pi x_1)\sin(2\pi x_2)$, with no explicit exact solution known in advance.
We therefore adopt the numerical solution (denoted by $U^n_{ref}$) obtained at the fine mesh with $\overline{h}=\frac{1}{200}, \tau=\frac{1}{400}$ as the reference solution, and the following formula to calculate the temporal rates at time $t_n$:
\[
\log_2\frac{\text{Error}(\tau;\overline{h},t_n)}{\text{Error}(\frac{\tau}{2};\overline{h},t_n)},\quad\text{where~~ }
\text{Error}(\tau;\overline{h},t_n)=\|U^n-U^n_{ref}\|.
\]
\par
Similar as in \textit{Example 2}, the errors and convergence rates are reported for $\alpha=0.3,0.6$ and $0.9$, respectively.
However, unlike the results in \textit{Example 2}, the SFTR-$\frac{1}{2}$ can still yield almost optimal convergence accuracy even through the solution is weak regular, which is much better than the classical F-BDF2 suffering a severe reduction in precision.
\begin{table}[]
\centering
\caption{Comparisons of our scheme (\ref{com.4}) using SFTR-$\frac{1}{2}$ with F-BDF2 for \textit{Example 3}.}\label{tab2}
{\footnotesize
{\renewcommand{\arraystretch}{1.4}
\begin{tabular*}{\hsize}{@{}@{\extracolsep{\fill}}cllllll@{}}
\toprule
\multirow{2}{*}{$\alpha$} & \multicolumn{1}{c}{\multirow{2}{*}{$\tau$}} & \multicolumn{2}{c}{SFTR-$\frac{1}{2}$}                &  & \multicolumn{2}{c}{F-BDF2}                            \\ \cline{3-4} \cline{6-7}
                          & \multicolumn{1}{c}{}                        & \multicolumn{1}{c}{Errer} & \multicolumn{1}{c}{Rates} &  & \multicolumn{1}{c}{Errer} & \multicolumn{1}{c}{Rates} \\ \hline
\multirow{4}{*}{0.3}      & 1/20                                        & 1.7143E-04                &                           &  & 4.0081E-02                &                           \\
                          & 1/40                                        & 5.7584E-05                & 1.57                      &  & 2.7325E-02                & 0.55                      \\
                          & 1/80                                        & 1.8625E-05                & 1.63                      &  & 1.7007E-02                & 0.68                      \\
                          & 1/160                                       & 5.1691E-06                & 1.85                      &  & 8.6550E-03                & 0.97                      \\ \hline
\multirow{4}{*}{0.6}      & 1/20                                        & 7.7902E-05                &                           &  & 1.5362E-02                &                           \\
                          & 1/40                                        & 1.8377E-05                & 2.08                      &  & 9.1160E-03                & 0.75                      \\
                          & 1/80                                        & 4.2841E-06                & 2.10                      &  & 4.9807E-03                & 0.87                      \\
                          & 1/160                                       & 8.8900E-07                & 2.27                      &  & 2.2460E-03                & 1.15                      \\ \hline
\multirow{4}{*}{0.9}      & 1/20                                        & 2.7149E-05                &                           &  & 4.6119E-03                &                           \\
                          & 1/40                                        & 6.6695E-06                & 2.03                      &  & 2.3217E-03                & 0.99                      \\
                          & 1/80                                        & 1.6462E-06                & 2.02                      &  & 1.0908E-03                & 1.09                      \\
                          & 1/160                                       & 3.5626E-07                & 2.21                      &  & 4.2951E-04                & 1.34                      \\
\bottomrule
\end{tabular*}}}
\end{table}
\subsection{Comparison with the L2-1$_\sigma$ formula}
\textit{Example 4.}\quad
In this example, we compare our scheme with that in \cite{LiaoTangZhou} formulated by the well-known L2-1$_\sigma$ formula on nonuniform meshes.
For better demonstration, we calculate the $L^2$ norm errors at the finial time $t=T$, under the condition that the number of grid points among $[0,T]$, which is denoted by $N$, is exactly the same for both schemes.
\par
\begin{figure}[htbp]
\centering
\subfigure[]{
\begin{minipage}[t]{0.45\linewidth}
\centering
\includegraphics[width=1\textwidth]{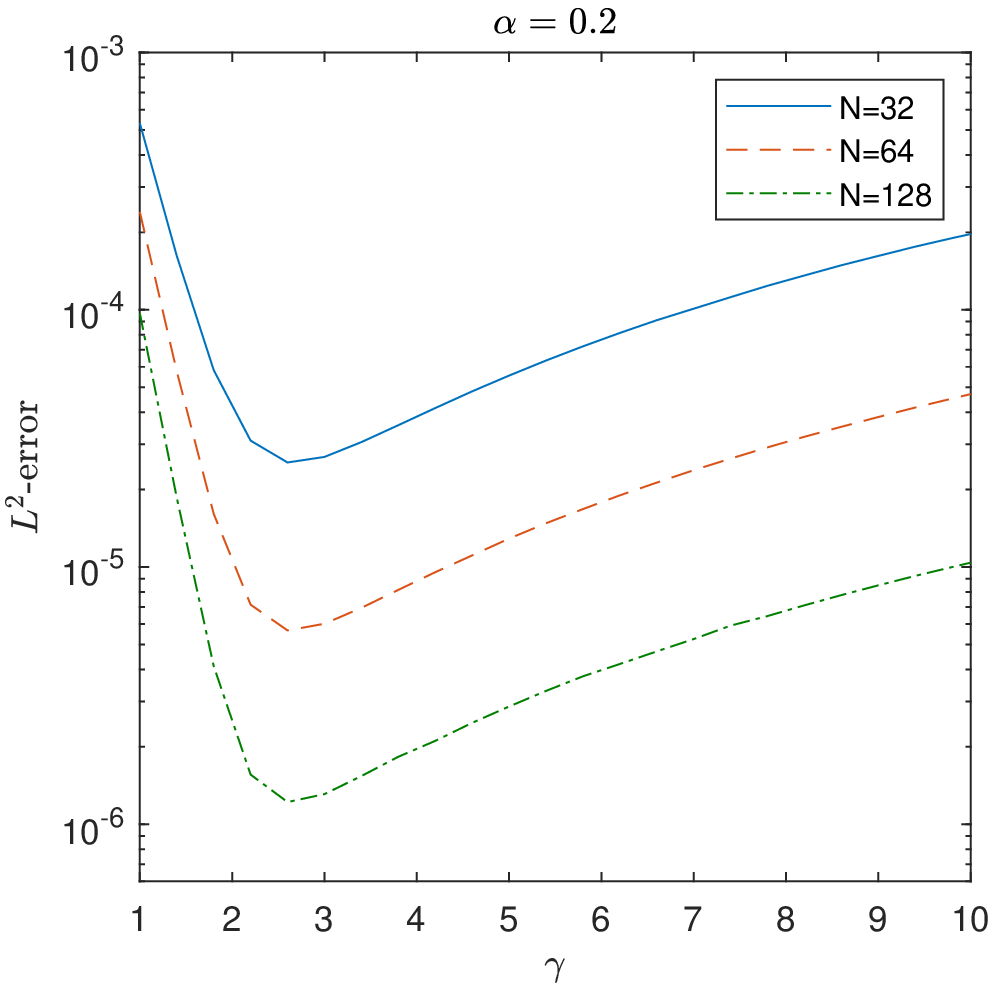}
\end{minipage}%
}%
\subfigure[]{
\begin{minipage}[t]{0.45\linewidth}
\centering
\includegraphics[width=1\textwidth]{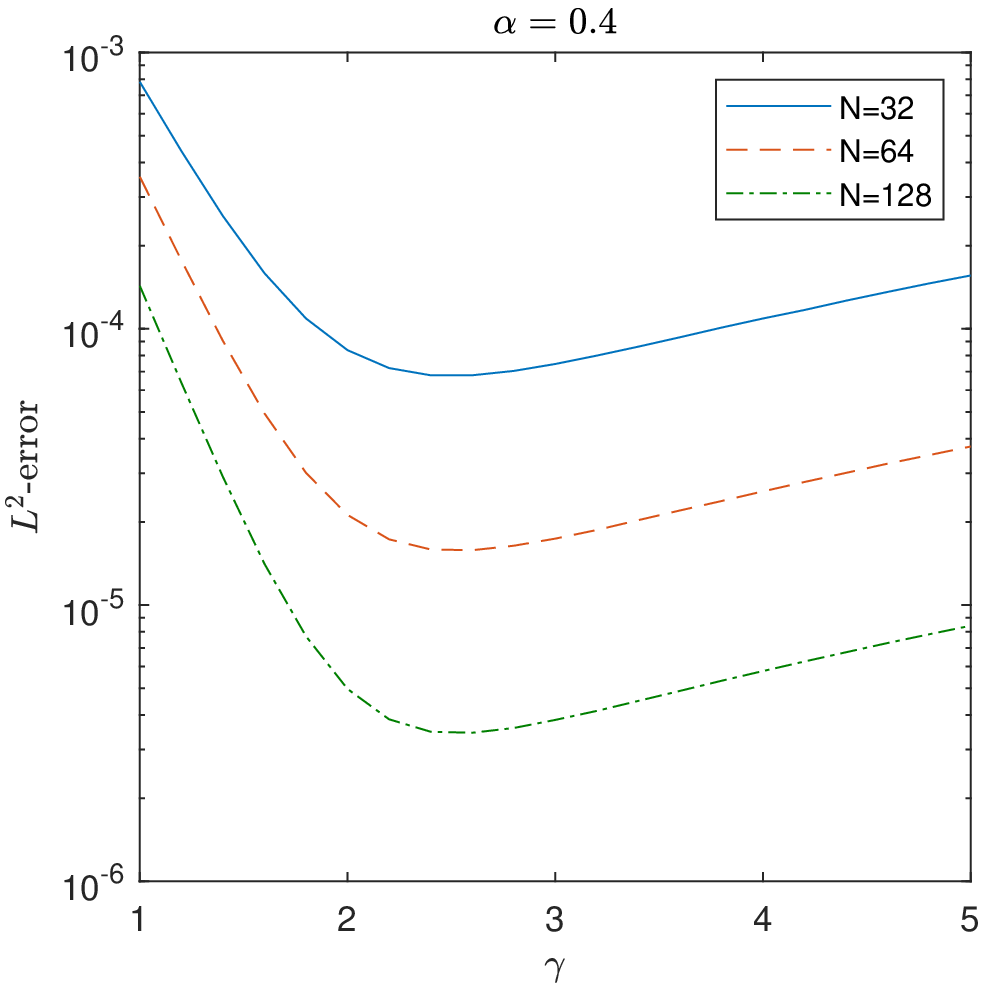}
\end{minipage}%
}%
\\
\centering
\subfigure[]{
\begin{minipage}[t]{0.45\linewidth}
\centering
\includegraphics[width=1\textwidth]{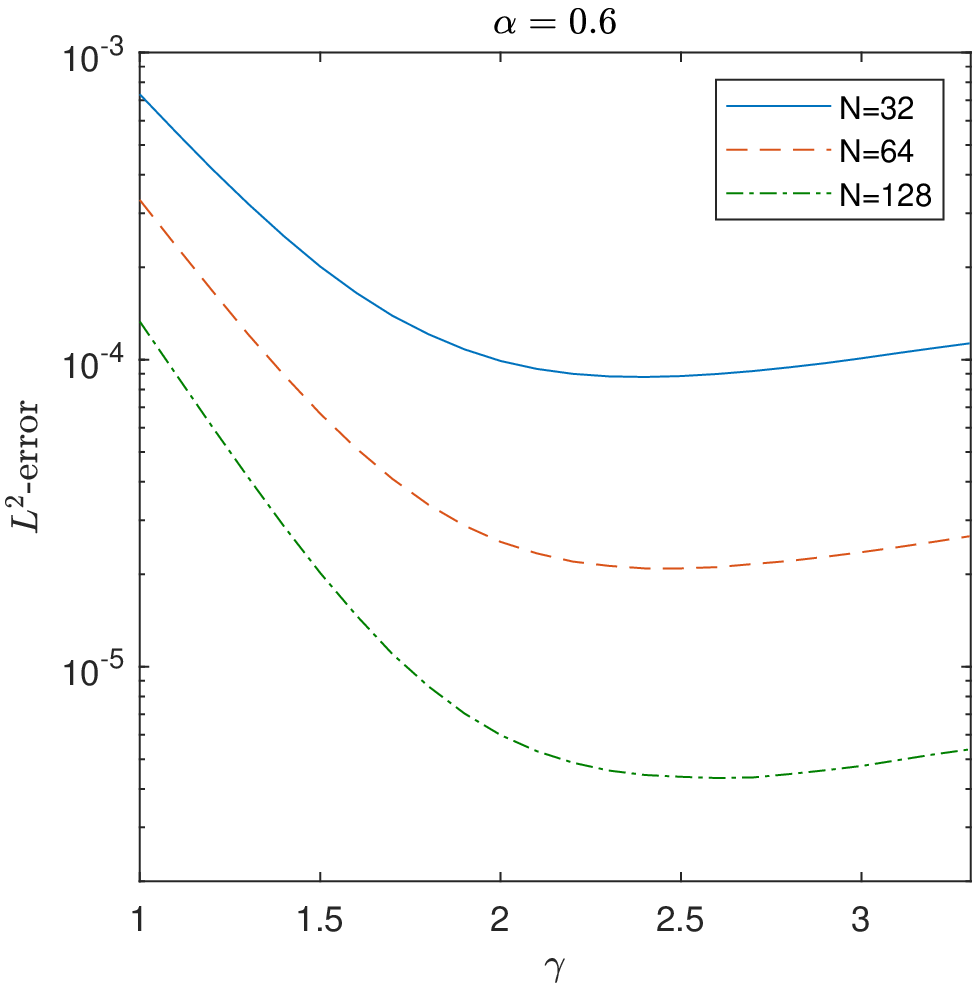}
\end{minipage}%
}%
\subfigure[]{
\begin{minipage}[t]{0.45\linewidth}
\centering
\includegraphics[width=1\textwidth]{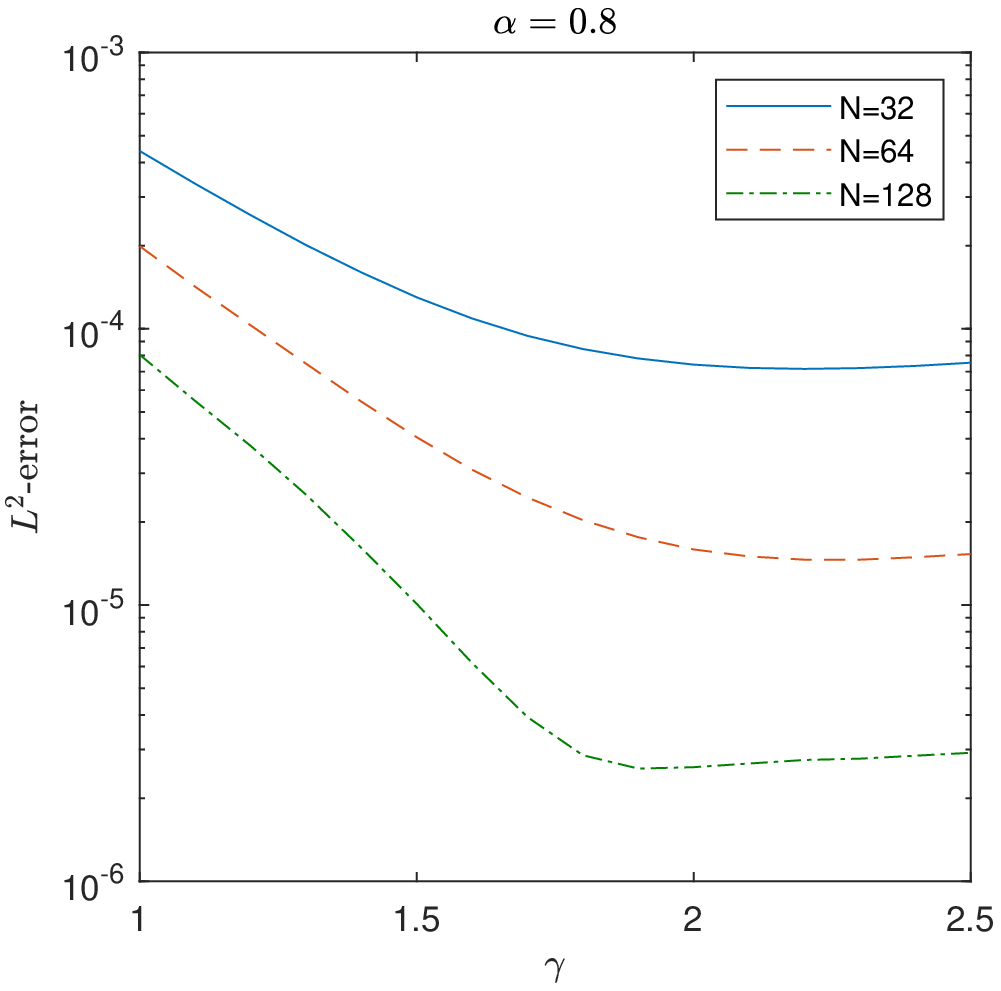}
\end{minipage}%
}%
\centering
\caption{Errors in $L^2$ norm of the scheme in \cite{LiaoTangZhou} for varying $\gamma \in [1,\frac{2}{\alpha}]$}\label{C3} with (a) $\alpha=0.2$, (b) $\alpha=0.4$, (c) $\alpha=0.6$ and (d) $\alpha=0.8$, respectively.
\end{figure}

Let $\Omega=(0,1)\times(0,1)$, $T=1$ and $\varepsilon=0.01$.
Choose the same initial condition as that in \textit{Example 3}.
By graded meshes we mean that for some given $\gamma\geq 1$, the interval $[0,T]$ is divided by points $t_n=T(n/N)^\gamma$.
Similar as \textit{Example 3}, the reference solutions $U_{ref}$, which depend generally on $\alpha,h,N$ and on $\gamma$ in addition if graded meshes are considered, are obtained under the fine mesh with $h=\overline{h}=\frac{1}{200}, N=400$.
Fig.\ref{C3} illustrates the errors in $L^2$ norm of the scheme \cite{LiaoTangZhou} for different $\gamma$ with $\gamma \in [1,\frac{2}{\alpha}]$.
Note that if $\gamma=1$, the mesh is uniform and if $\gamma =\frac{2}{\alpha}$, the uniform second-order convergence rate is arrived at according to \cite{LiaoTangZhou}.
One easily observes from Fig.\ref{C3} that to obtain the minimum error at final time $t=T$, the mesh parameter $\gamma$ should be between $1$ and $\frac{2}{\alpha}$, e.g., from Fig.\ref{C3} (a) where $\alpha=0.2$, $\gamma=2.6$ results in the minimum error.
We therefore in Table \ref{tab3} collect some key results and denote by ``L2-1$_\sigma$ (uniform)'' the errors on uniform mesh (i.e., $\gamma=1$) and by ``L2-1$_\sigma$ (min)'' the minimum errors among $\gamma \in [1,\frac{2}{\alpha}]$.
Obviously, SFTR-$\frac{1}{2}$ is much more accurate than L2-1$_\sigma$ on uniform meshes and both of these two methods are evenly matched (so long as $\alpha$ is not too small) if graded meshes are equipped to L2-1$_\sigma$.
We also remark that SFTR-$\frac{1}{2}$ performs even better for large $\alpha$ (e.g., $\alpha \geq 0.4$ in our example), though L2-1$_\sigma$ is more accurate when $\alpha$ is small in which case more work must be involved to avoid the round-off errors polluting the coefficients of the formula (see (3.2) and (3.3) in \cite{LiaoTangZhou}).
\begin{table}[]
\centering
\caption{Comparisons of SFTR-$\frac{1}{2}$ with L2-1$_\sigma$ for \textit{Example 4}.}\label{tab3}
{\footnotesize
{\renewcommand{\arraystretch}{1.4}
\begin{tabular*}{\hsize}{@{}@{\extracolsep{\fill}}llllllll@{}}
\toprule
\multicolumn{1}{c}{}                   & \multicolumn{3}{c}{$\alpha=0.2$} & \multicolumn{1}{c}{} & \multicolumn{3}{c}{$\alpha=0.4$} \\ \cline{2-4} \cline{6-8}
                                       & $N=32$    & $N=64$    & $N=128$  &                      & $N=32$    & $N=64$    & $N=128$  \\ \hline
SFTR-$\frac{1}{2}$                     & 7.00E-05  & 2.73E-05  & 9.01E-06 &                      & 3.14E-05  & 9.11E-06  & 2.48E-06 \\
L2-1$_\sigma$(uniform)               & 5.31E-04  & 2.40E-04  & 9.68E-05 &                      & 7.83E-04  & 3.55E-04  & 1.43E-04 \\
L2-1$_\sigma$($\min$) & 2.55E-05  & 5.67E-06  & 1.22E-06 &                      & 6.79E-05  & 1.58E-05  & 3.45E-06 \\ \hline
                                       & \multicolumn{3}{c}{$\alpha=0.6$} &                      & \multicolumn{3}{c}{$\alpha=0.8$} \\ \cline{2-4} \cline{6-8}
                                       & $N=32$    & $N=64$    & $N=128$  &                      & $N=32$    & $N=64$    & $N=128$  \\ \hline
SFTR-$\frac{1}{2}$                     & 1.37E-05  & 3.42E-06  & 7.70E-07 &                      & 8.66E-06  & 2.10E-06  & 5.28E-07 \\
L2-1$_\sigma$(uniform)               & 7.31E-04  & 3.31E-04  & 1.33E-04 &                      & 4.40E-04  & 1.99E-04  & 8.05E-05 \\
L2-1$_\sigma$($\min$) & 8.78E-05  & 2.09E-05  & 4.34E-06 &                      & 7.16E-05  & 1.46E-05  & 2.56E-06 \\ \bottomrule
\end{tabular*}}}
\end{table}
\section{Conclusion}\label{sec.conc}
{In this work, we develop a finite difference scheme for time-fractional Allen-Cahn equation by SFTR-$\frac{1}{2}$ for time discretization and the centered difference formula for space variables.
We first prove some key properties of the weights of SFTR-$\frac{1}{2}$ and further show that the scheme proposed is uniquely solvable.
Moreover, the discrete energy decay property of the proposed scheme is exhibited by using some properties of two sets of coefficients $\theta_k$ and $\vartheta_k$, which are related to $\omega_k$.
The discrete maximum-principle preserved of the scheme is also proved rigorously, by resorting to the negativity of the weights $\omega_k$ for $k\geq 1$.
Numerical results confirm the correctness of our theoretical results concerning discrete energy decay law and maximum-principle, and reveal that SFTR-$\frac{1}{2}$ can resolve automatically the intrinsic singularity of tFAC equation, indicating SFTR-$\frac{1}{2}$ is superior to other high-order approximation formulas on uniform meshes (e.g., the fractional BDF2 or L2-1$_\sigma$).
Comparisons with L2-1$_\sigma$ method under graded meshes also confirm the high accuracy of our method.
\par
There are some aspects deserving further study, e.g., proposing the theoretical analysis for the optimal convergence results of SFTR-$\frac{1}{2}$ when applied to nonlinear fractional equations, which will be our future work. }

\section*{Acknowledgments} 
The authors are grateful to two anonymous referees and editors for their valuable suggestions which improve the presentation of this work greatly.
This work is supported by the Natural Science Foundation of Inner Mongolia Autonomous Region of China (No. 2021BS01003 to G.Z.), NSF of China (Nos. 12171177 and 12011530058 to C.H., 12201322 to B.Y.), RFBR (No. 20-51-53007 to A. A.) and North-Caucasus Center for Mathematical Research under agreement No. 075-02-2021-1749 with the Ministry of Science and Higher Education of the Russian Federation.

 \section*{Conflict of interest}

 The authors declare that they have no conflict of interest.

\bibliographystyle{spmpsci}      
\bibliography{reference}   

\begin{thebibliography}{10}
\providecommand{\url}[1]{{#1}}
\providecommand{\urlprefix}{URL }
\expandafter\ifx\csname urlstyle\endcsname\relax
  \providecommand{\doi}[1]{DOI~\discretionary{}{}{}#1}\else
  \providecommand{\doi}{DOI~\discretionary{}{}{}\begingroup
  \urlstyle{rm}\Url}\fi

\bibitem{MaskariKaraa}
Al-Maskari, M., Karaa, S.: {The time-fractional Cahn--Hilliard equation:
  analysis and approximation}.
\newblock IMA J. Numer. Anal.  (2021)

\bibitem{alikhanov2015new}
Alikhanov, A.A.: {A new difference scheme for the time fractional diffusion
  equation}.
\newblock J. Comput. Phys. \textbf{280}, 424--438 (2015)

\bibitem{AllenCahn}
Allen, S.M., Cahn, J.W.: {A microscopic theory for antiphase boundary motion
  and its application to antiphase domain coarsening}.
\newblock Acta metall. \textbf{27}(6), 1085--1095 (1979)

\bibitem{AndersonMcFaddenWheeler}
Anderson, D.M., McFadden, G.B., Wheeler, A.A.: {Diffuse-interface methods in
  fluid mechanics}.
\newblock Annu. Rev. Fluid Mech. \textbf{30}(1), 139--165 (1998)

\bibitem{bu2017finite}
Bu, W., Xiao, A., Zeng, W.: {Finite difference/finite element methods for
  distributed--order time fractional diffusion equations}.
\newblock J. Sci. Comput. \textbf{72}(1), 422--441 (2017)

\bibitem{ChenZhangZhaoCaoWangZhang}
Chen, L., Zhang, J., Zhao, J., Cao, W., Wang, H., Zhang, J.: {An accurate and
  efficient algorithm for the time-fractional molecular beam epitaxy model with
  slope selection}.
\newblock Comput. Phys. Commun. \textbf{245}, 106842 (2019)

\bibitem{chen2020spectrally}
Chen, S., Shen, J., Zhang, Z., Zhou, Z.: {A spectrally accurate approximation
  to subdiffusion equations using the log orthogonal functions}.
\newblock SIAM J. Sci. Comput. \textbf{42}(2), A849--A877 (2020)

\bibitem{DuYangZhou}
Du, Q., Yang, J., Zhou, Z.: {Time-fractional Allen--Cahn equations: analysis
  and numerical methods}.
\newblock J. Sci. Comput. \textbf{85}(2), 1--30 (2020)

\bibitem{gao2011compact}
Gao, G., Sun, Z.: {A compact finite difference scheme for the fractional
  sub-diffusion equations}.
\newblock J. Comput. Phys. \textbf{230}(3), 586--595 (2011)

\bibitem{hou2021highly}
Hou, D., Xu, C.: {Highly efficient and energy dissipative schemes for the time
  fractional Allen-Cahn equation}.
\newblock SIAM J. Sci. Comput. \textbf{43}(5), A3305--A3327 (2021)

\bibitem{hou2022second}
Hou, D., Xu, C.: {A second order energy dissipative scheme for time fractional
  $L^2$ gradient flows using SAV approach}.
\newblock J. Sci. Comput. \textbf{90}(1), 1--22 (2022)

\bibitem{HouZhuXu}
Hou, D., Zhu, H., Xu, C.: {Highly efficient schemes for time-fractional
  Allen-Cahn equation using extended SAV approach}.
\newblock Numer. Algorithms pp. 1--32 (2021)

\bibitem{HouTang01}
Hou, T., Tang, T., Yang, J.: {Numerical analysis of fully discretized
  Crank--Nicolson scheme for fractional-in-space Allen--Cahn equations}.
\newblock J. Sci. Comput. \textbf{72}(3), 1214--1231 (2017)

\bibitem{huang2020optimal}
Huang, C., Stynes, M.: {Optimal $H^1$ spatial convergence of a fully discrete
  finite element method for the time-fractional Allen-Cahn equation}.
\newblock Adv. Comput. Math. \textbf{46}(4), 1--20 (2020)

\bibitem{ji2020adaptive}
Ji, B., Liao, H., Gong, Y., Zhang, L.: {Adaptive linear second-order energy
  stable schemes for time-fractional Allen-Cahn equation with volume
  constraint}.
\newblock Commun. Nonlinear Sci. Numer. Simul. \textbf{90}, 105366 (2020)

\bibitem{JiaZhangXuJiang}
Jia, J., Zhang, H., Xu, H., Jiang, X.: {An efficient second order stabilized
  scheme for the two dimensional time fractional Allen-Cahn equation}.
\newblock Appl. Numer. Math. \textbf{165}, 216--231 (2021)

\bibitem{jiang2017fast}
Jiang, S., Zhang, J., Zhang, Q., Zhang, Z.: {Fast evaluation of the Caputo
  fractional derivative and its applications to fractional diffusion
  equations}.
\newblock Commun. Comput. Phys. \textbf{21}(3), 650--678 (2017)

\bibitem{JinLazarovZhou}
Jin, B., Lazarov, R., Zhou, Z.: {Numerical methods for time-fractional
  evolution equations with nonsmooth data: a concise overview}.
\newblock Comput. Methods Appl. Mech. Engrg. \textbf{346}, 332--358 (2019)

\bibitem{jin2018numerical}
Jin, B., Li, B., Zhou, Z.: {Numerical analysis of nonlinear subdiffusion
  equations}.
\newblock SIAM J. Numer. Anal. \textbf{56}(1), 1--23 (2018)

\bibitem{kopteva2019error}
Kopteva, N.: {Error analysis of the L1 method on graded and uniform meshes for
  a fractional--derivative problem in two and three dimensions}.
\newblock Math. Comp. \textbf{88}(319), 2135--2155 (2019)

\bibitem{li2018finite}
Li, C., Yi, Q.: {Finite difference method for two-dimensional nonlinear
  time-fractional subdiffusion equation}.
\newblock Fract. Calc. Appl. Anal. \textbf{21}(4), 1046--1072 (2018)

\bibitem{li2021novel}
Li, D., Sun, W., Wu, C.: {A novel numerical approach to time-fractional
  parabolic equations with nonsmooth solutions}.
\newblock Numer. Math. Theory Methods Appl. \textbf{14}(2), 355--376 (2021)

\bibitem{li2018unconditionally}
Li, D., Zhang, J., Zhang, Z.: {Unconditionally optimal error estimates of a
  linearized Galerkin method for nonlinear time fractional
  reaction--subdiffusion equations}.
\newblock J. Sci. Comput. \textbf{76}(2), 848--866 (2018)

\bibitem{li2019nonconforming}
Li, M., Zhao, J., Huang, C., Chen, S.: {Nonconforming virtual element method
  for the time fractional reaction--subdiffusion equation with non-smooth
  data}.
\newblock J. Sci. Comput. \textbf{81}(3), 1823--1859 (2019)

\bibitem{li2021conforming}
Li, M., Zhao, J., Huang, C., Chen, S.: {Conforming and nonconforming VEMs for
  the fourth-order reaction--subdiffusion equation: a unified framework}.
\newblock IMA J. Numer. Anal.  (2021)

\bibitem{LiaoTangZhou}
Liao, H., Tang, T., Zhou, T.: {A second-order and nonuniform time-stepping
  maximum-principle preserving scheme for time-fractional Allen-Cahn
  equations}.
\newblock J. Comput. Phys. \textbf{414}, 109473 (2020)

\bibitem{LiaoTangZhou1}
Liao, H., Tang, T., Zhou, T.: {An energy stable and maximum bound preserving
  scheme with variable time steps for time fractional Allen--Cahn equation}.
\newblock SIAM J. Sci. Comput. \textbf{43}(5), A3503--A3526 (2021)

\bibitem{LiaoZhuWang}
Liao, H., Zhu, X., Wang, J.: {The variable-step L1 scheme preserving a
  compatible energy law for time-fractional Allen-Cahn equation}.
\newblock arXiv preprint arXiv:2102.07577  (2021)

\bibitem{liu2009numerical}
Liu, F., Yang, C., Burrage, K.: {Numerical method and analytical technique of
  the modified anomalous subdiffusion equation with a nonlinear source term}.
\newblock J. Comput. Appl. Math. \textbf{231}(1), 160--176 (2009)

\bibitem{LiuChengWang}
Liu, H., Cheng, A., Wang, H., Zhao, J.: {Time-fractional Allen--Cahn and
  Cahn--Hilliard phase-field models and their numerical investigation}.
\newblock Comput. Math. Appl. \textbf{76}(8), 1876--1892 (2018)

\bibitem{YinLiuLiZhang444}
Liu, Y., Yin, B., Li, H., Zhang, Z.: {The unified theory of shifted convolution
  quadrature for fractional calculus}.
\newblock J. Sci. Comput. \textbf{89}(1), 1--24 (2021)

\bibitem{Lubich1}
Lubich, C.: {Discretized fractional calculus}.
\newblock SIAM J. Math. Anal. \textbf{17}(3), 704--719 (1986)

\bibitem{mohebbi2013high}
Mohebbi, A., Abbaszadeh, M., Dehghan, M.: {A high-order and unconditionally
  stable scheme for the modified anomalous fractional sub-diffusion equation
  with a nonlinear source term}.
\newblock J. Comput. Phys. \textbf{240}, 36--48 (2013)

\bibitem{QuanTangYang}
Quan, C., Tang, T., Yang, J.: {How to define dissipation-preserving energy for
  time-fractional phase-field equations}.
\newblock arXiv preprint arXiv:2007.14855  (2020)

\bibitem{ShaoRappelLevine}
Shao, D., Rappel, W.J., Levine, H.: {Computational model for cell
  morphodynamics}.
\newblock Phys. Rev. Lett. \textbf{105}(10), 108104 (2010)

\bibitem{ShenTangYang}
Shen, J., Tang, T., Yang, J.: {On the maximum principle preserving schemes for
  the generalized Allen--Cahn equation}.
\newblock Commun. Math. Sci. \textbf{14}(6), 1517--1534 (2016)

\bibitem{Stynes}
Stynes, M.: {Too much regularity may force too much uniqueness}.
\newblock Fract. Calc. Appl. Anal. \textbf{19}(6), 1554--1562 (2016)

\bibitem{stynes2017error}
Stynes, M., O'Riordan, E., Gracia, J.L.: {Error analysis of a finite difference
  method on graded meshes for a time-fractional diffusion equation}.
\newblock SIAM J. Numer. Anal. \textbf{55}(2), 1057--1079 (2017)

\bibitem{TangYuZhou}
Tang, T., Yu, H., Zhou, T.: On energy dissipation theory and numerical
  stability for time-fractional phase-field equations.
\newblock SIAM J. Sci. Comput. \textbf{41}(6), A3757--A3778 (2019)

\bibitem{XuPrinceSnakes}
Xu, C., Prince, J.L.: {Snakes, shapes, and gradient vector flow}.
\newblock IEEE Trans. Image Process. \textbf{7}(3), 359--369 (1998)

\bibitem{Yin01}
Yin, B., Liu, Y., Li, H.: {Necessity of introducing non-integer shifted
  parameters by constructing high accuracy finite difference algorithms for a
  two-sided space-fractional advection--diffusion model}.
\newblock Appl. Math. Lett. \textbf{105}, 106347 (2020)

\bibitem{YinLiuLiZhang333}
Yin, B., Liu, Y., Li, H., Zhang, Z.: {Two families of novel second-order
  fractional numerical formulas and their applications to fractional
  differential equations}.
\newblock arXiv preprint arXiv:1906.01242  (2019)

\bibitem{YinLiuLiZhang222}
Yin, B., Liu, Y., Li, H., Zhang, Z.: {Efficient shifted fractional trapezoidal
  rule for subdiffusion problems with nonsmooth solutions on uniform meshes}.
\newblock BIT pp. 1--36 (2021)

\bibitem{zhang2020efficient}
Zhang, H., Yang, X., Xu, D.: {An efficient spline collocation method for a
  nonlinear fourth-order reaction subdiffusion equation}.
\newblock J. Sci. Comput. \textbf{85}(1), 1--18 (2020)

\bibitem{zhao2014collocation}
Zhao, J., Xiao, J., Ford, N.J.: {Collocation methods for fractional
  integro-differential equations with weakly singular kernels}.
\newblock Numer. Algorithms \textbf{65}(4), 723--743 (2014)

\end{thebibliography}


\end{document}